\pdfoutput=1
\documentclass[preprin]{elsarticle}
\usepackage{Styles/ART-DFN-Adapt}
\usepackage{xcolor}
\usepackage{listings}

{\left\lbrace\begin{array}{@{}l@{}}}%
{\end{array}\right.}

\definecolor{CommentGreen}{rgb}{0.16,0.55,0.20}
\definecolor{WordBlue}{rgb}{0.08,0.28,0.6}
\definecolor{DefBlue}{rgb}{0.09,0.31,0.8}

\graphicspath{{Immagini/}}

\begin{document}
\begin{frontmatter}
  \title{Refinement strategies for polygonal meshes applied to adaptive VEM discretization.\tnoteref{t1}}%
  \tnotetext[t1]{This research has been supported by INdAM-GNCS
    Projects 2018 and 2019, by the MIUR project PRIN 201744KLJL\_004,
    and by the MIUR project ``Dipartimenti di Eccellenza 2018-2022''
    (CUP E11G18000350001).  Computational resources were partially
    provided by HPC@POLITO (\url{http://hpc.polito.it}).}%

  \author[poli,indam]{Stefano Berrone\corref{cor1}}%
  \ead{stefano.berrone@polito.it}%
  \author[poli,indam]{Andrea Borio}%
  \ead{andrea.borio@polito.it}%
  \author[poli,indam]{Alessandro D'Auria}%
  \ead{alessandro.dauria@polito.it}%
  \address[poli]{Dipartimento di Scienze Matematiche, Politecnico di
    Torino\\Corso Duca degli Abruzzi 24, Torino, 10129, Italy}%
  \address[indam]{Member of the INdAM research group GNCS}%
  \cortext[cor1]{Corresponding author}%

  \begin{abstract}
    In the discretization of differential problems on complex
    geometrical domains, discretization methods based on polygonal and
    polyhedral elements are powerful tools. Adaptive mesh refinement
    for such kind of problems is very useful as well and states new issues, here tackled, concerning good quality mesh elements and reliability of the simulations.
    In this paper we numerically investigate optimality with respect
    to the number of degrees of freedom of the numerical solutions
    obtained by the different refinement strategies proposed.
    A geometrically
    complex geophysical problem is used as test problem for several
    general purpose and problem dependent refinement strategies.
  \end{abstract}

  \begin{keyword}
    Mesh adaptivity \sep Polygonal mesh refinement \sep Virtual
    Element Method \sep Discrete Fracture Network flow simulations
    \sep Simulations in complex geometries \sep A posteriori error
    estimates

    \MSC[2010] 65N30 \sep 65N50 \sep 68U20 \sep 86-08 \sep 86A05
  \end{keyword}

\end{frontmatter}

\section{Introduction}
\label{sec:intro}

In last years, a growing interest has arised for the development of
numerical methods for the solution of partial differential equations
using general polygonal meshes. These methods are well suited
for handling domains featuring geometrical complexities that
can yield situations where the generation of good quality conforming meshes
can be particularly expensive or even unfeasible.
On the other hand, the use of
polygonal and polyhedral elements in conjunction with adaptive mesh
refinement states new problems. In this paper we investigate several refinement
strategies for polygonal meshes.

As model problem to test the different refinement strategies
we consider the simulation of flow in a fractured medium.
Flow simulation in underground fractured media is a key aspect
in many applications such as aquifers monitoring, nuclear waste
disposal, $CO_2$ geological storage or control of contaminant dispersion
in the subsoil.


In those situations where the rock matrix can be considered perfectly
impervious, the Discrete Fracture Network (DFN) model is often used,
representing fractures as planar polygons randomly intersected  in the
3D space \cite{DE88,Fidelibus2009,CLMTBDFP,Garipov}.  These kind of
domains usually present many geometrical challenges, in particular
when a conforming mesh is required, 
for example in order to apply standard finite
element approaches. To circumvent these difficulties, many approaches
have been devised, relying on domain decomposition techniques
\cite{P1,P2}, or devising particular meshing strategies
\cite{karimi2014unstructured,Gable,FOURNO2019713,NGO201766,NOETINGER2015205},
or using extensions of the finite element method \cite{
  BBSmixed,FS13,formaggia2014optimal,SF2014,Antonietti2016,FKS,Lenti2003,DF99,al2017generalized,flemisch2018benchmarks,Masson2017,jaffre2012modeling,AHMED2017265}.
In \cite{BPSa,BPSb,BPSc}, a PDE-constrained optimization approach was
devised, enabling for the use of completely non-conforming meshes,
where finite element, extended finite element \cite{BFPS} or virtual
element spaces \cite{BBPS} can be used. A scalable parallel
implementation of such approach is discussed in \cite{BSV} and a GPU accelerated version in \cite{BDV}. In
\cite{BBoS} and \cite{BERRONE2019904}, a residual a posteriori error
estimate is devised and applied to large scale DFNs.  Here we follow
the approach based on the Virtual Element Method (VEM)
\cite{Beirao2013a,Ahmad2013,Beirao2015b,Beirao2017} first introduced
for DFN flow simulations in \cite{BBPS} and extended in
\cite{BBBPS,BBS,BBapost}.

The huge cost of large scale simulations on the scale of a geological
basin raises the issue of efficiency. Moreover, due to the large uncertainty
in geometrical configurations and hydrogeological parameters a
stochastic approach is advisable yielding to the need of
many simulations.
Beside efficiency, reliability of these
simulations is of crucial importance being connected with risk
analysis of human activities.
Efficiency and reliability strongly affect the
application of uncertainty quantification techniques for the
estimation of relevant quantities of interest \cite{BCPSa,BCPSb}.
For these reasons adaptivity plays a key role both on the side of
reliability and efficiency.

Although we focus on the simulation of flow in DFN using VEM, the
refinement strategies explored here are suitable for any kind of
method relying on polygonal meshes, except for the so-called ``Trace
Direction'' refinement, that is particularly suited for DFN problems taking into account the behaviour of the solution approaching the traces. The approach
here followed is based on isotropic ``a posteriori'' error estimates
and aims at generating good quality isotropic polygonal
elements. Anisotropic estimates and mesh refinement strategies are currently
under investigation in \cite{AntoniettiXXXX}.

In Section \ref{sec:dfn} we introduce the DFN model used to test the refinement strategies and define some useful notations. In Section \ref{sec:vem} we introduce its Virtual Element discretization, we define the error estimators and state the ``a posteriori'' error estimate. In Sections \ref{sec:minimal-mesh} and
\ref{sec:algorithm} we describe the initial mesh generation, the
adaptive algorithm, and the polygonal refinement strategies. In
Section \ref{sec:num_res_TP} we validate their application to some test
cases with known solution, and, in Section \ref{sec:num_res_DFN}, we
analyse some results on a more realistic DFN.

 
\section{Discrete Fracture Networks}
\label{sec:dfn}

Let us consider a set of open convex planar polygonal fractures
$F_i \subset \mathbb{R}^3$ with $i=1,...,N$, with boundary
$\partial F$. A DFN is $\Omega=\bigcup_{i}F_i$, with boundary
$\partial \Omega$. Even though the fractures are planar, their
orientations in space are arbitrary, such that $\Omega$ is a 3D
set. The set $\Gamma_D\subset \partial \Omega$ is where Dirichlet
boundary conditions are imposed, and we assume
$\Gamma_D\neq \emptyset$, whereas
$\Gamma_N=\partial \Omega \setminus \Gamma_D$, is the portion of the
boundary with Neumann boundary conditions. Dirichlet and Neumann
boundary conditions are prescribed by the functions
$h^D \in \sobh{\frac12}{\Gamma_{D}}$ and
$g^N\in\sobh{-\frac12}{\Gamma_{N}}$ on the Dirichlet and Neumann part
of the boundary, respectively. We further set
$\Gamma_{iD} = \Gamma_{D} \cap \partial F_i$,
$\Gamma_{iN} = \Gamma_{N} \cap \partial F_i$, and
$h_i^D={h^D}_{|\Gamma_{iD}}$ and $g_i^N={g^N}_{|\Gamma_{iN}}$. The set
$\CAL{S}$ collects all the traces, i.e. the intersections between
fractures, and each trace $S\in \CAL{S}$ is given by the intersection
of exactly two fractures, $ S=\bar{F}_i \cap \bar{F}_j$, such that
there is a one to one relationship between a trace $S$ and a couple of
fracture indices $\{i,j\}=\CAL{I}(S)$. We will also denote by
$\CAL{S}_i$ the set of traces belonging to fracture $F_i$.

Subsurface flow is governed by the gradient of the hydraulic head
$H=\CAL{P}+\zeta$, where $\CAL{P}=p/(\varrho g)$ is the pressure head,
$p$ is the fluid pressure, $g$ is the gravitational acceleration
constant, $\varrho$ is the fluid density and $\zeta$ is the elevation.

We define the following functional spaces:
\begin{align*}
  V_i&=\sobhob{1}{\rm 0}{F_i}=
       \left\{
       v\in \sobh{1}{F_i}: v_{|_{\Gamma_{iD}}}=0
       \right\},
  \\
  V_i^D&=\sobhob{1}{\rm D}{F_i}=
         \left\{
         v\in \sobh{1}{F_i}: v_{|_{\Gamma_{iD}}}=h_{i}^D
         \right\},
\end{align*}
and
\begin{equation*}
    V = \left\{ v\colon \eval{v}{F_i} \in V_i, \ \forall
      i\!=\!1,\ldots,N, \gamma_S (v_{|F_i})\!=\!\gamma_S (v_{|F_j}), \; \forall
      S\!\in\! \CAL{S}_i, \ \{i,j\}\!=\!\CAL{I}(S) \right\} \,,
\end{equation*}
\begin{equation*}
    V^D =\left\{ v\colon \eval{v}{F_i} \in V_i^D\!\!, \ \forall
      i\!=\!1,\ldots,N, \gamma_S (v_{|F_i})\!=\!\gamma_S (v_{|F_j}), \; \forall
      S\!\in\! \CAL{S}_i, \ \{i,j\}\!=\!\CAL{I}(S) \right\} \,,
\end{equation*}
where $\gamma_S$ is the trace operator onto $S$.
In order to formulate the DFN flow problem, let us define
$\a[i]{}{}\colon V_i^D \times V_i \to \mathbb{R}$ be defined as
\begin{equation*}
  \a[i]{w}{v} = \scal[F_i]{\K_i \nabla
    \left( \eval{w}{F_i} \right)}{\nabla\left( \eval{v}{F_i} \right)}
  \quad \forall w\in V^D,\ v \in V\,,
\end{equation*}
where $\mathcal{K}_i$ is the fracture transmissivity tensor, that we
assume to be constant on each fracture.
Let us denote by $H\in V^D$ the hydraulic head on the DFN and by $H_i\in V_i^D$ its restriction to the fracture $F_i$ and by
$v\in V$ the test functions.
Assuming the hydraulic head modelled by the Darcy law, the whole problem on the DFN is: find 
$H\in V^D$ such that $\forall v\in V$
\begin{equation}
  \label{eq:darcy-model}
  \sum_{i=1}^N \a[i]{H_i}{v} = \sum_{i=1}^N  \scal[F_i]{f_i }{v}+
                               \dual[\sobh{-\frac12}{\Gamma_{N_i}},
                               \sobh{\frac12}{\Gamma_{N_i}}]
                               {g_{i}^{N}}{\eval{v}{\Gamma_{N_i}}}.                 
\end{equation}
 
\section{Virtual Element discretization}
\label{sec:vem}

In this section we describe the Virtual Element discretization of
equation \eqref{eq:darcy-model} assuming a globally conforming
polygonal mesh is given on the DFN. A globally conforming mesh ia a polygonal mesh on each fractures such that polygon edges match exactly at traces. In
Section \ref{sec:minimal-mesh} we devise a way to obtain a mesh of
this type and Section \ref{sec:algorithm} describes some refinement
strategies that preserve the global conformity of the mesh.

Let $\Th$ be a globally conforming polygonal mesh of $\Omega$
fulfilling the regularity requirements needed by the Virtual Element
method \cite{Beirao2015b}, $E\in\Th$
be any polygon of this tessellation. Let $\Poly{k}{E}$ be the space of
polynomials of degree $\leq k$ defined on $E$. To define the discrete
functional space on $E$, we introduce the $\sobh{1}{}$-orthogonal
projector $\proj[\nabla]{k,E}{}\colon \sobh{1}{E} \to \Poly{k}{E}$
such that
\begin{align*}
  \scal[E]{\nabla\left(v-\proj[\nabla]{k,E}{} v\right)}{\nabla p}=0 \,, \forall
  p\in\Poly{k}{E}
  \\
  \intertext{and}
  \begin{cases}
    \scal[\partial E]{\proj[\nabla]{k,E}{} v}{1} = \scal[\partial
    E]{v}{1} & \text{if $k=1$}\,,
    \\
    \scal[E]{\proj[\nabla]{k,E}{} v}{1} = \scal[E]{v}{1} & \text{if
      $k\geq 1$}\,.
  \end{cases}
\end{align*}
The Virtual Element space of order $k\in\mathbb{N}$ on $E$ is defined
as
\begin{multline*}
  \label{eq:defVhE}
  \Vh^E = \left\{ v\in \sobh{1}{E}\colon\Delta
    v\in\Poly{k}{E},\,v\in\Poly{k}{e} \,\forall e\subset\partial E,
    \trace{v}{\partial E}{} \in \cont{\partial E}\,, \right.
  \\
  \quad\left. \scal[E]{v}{p}=\scal[E]{\proj[\nabla]{k,E}{}
      v}{p}\,\forall p\in\Poly{k}{E}/ \Poly{k-2}{E} \right\} \,,
\end{multline*}
where $\Poly{k}{E}/ \Poly{k-2}{E}$ denotes the subspace of
$\Poly{k}{E}$ containing polynomials that are $\lebl{E}$-orthogonal to
$\Poly{k-2}{E}$. Furthermore, let us denote by $\Th[i]$ the
restriction of $\Th$ to fracture $F_i$. The Virtual Element space on
$\Th[i]$ is
\begin{equation*}
  \Vh[i] = \left\{ v\in \cont{F_i} \colon v\in \Vh^E \;
    \forall E\in\Th[i] \right\} \,.
\end{equation*}
Since $\Th$ is globally conforming, we can define the global
discrete spaces as
\begin{align*}
  \Vh &= \left\{ v \in V \colon v_{|F_i} \in \Vh[i] \right\} \,,
  &
    \Vh^D &= \left\{ v \in V^D \colon v_{|F_i} \in \Vh[i] \right\} \,.
\end{align*}
A function in the above space is uniquely identified by the following
set of degrees of freedom:
\begin{enumerate}
\item the values at the vertices of the polygons;
\item if $k\geq 2$, for each edge $e$ of the mesh, the value of $v$ at
  $k-1$ internal points of $e$;
\item if $k\geq 2$, for each $E\in\Th$, the scaled moments
  $\frac{1}{\abs{E}}\scal[E]{v}{m_{\boldsymbol{\alpha}}}$ for all the
  scaled monomials $m_{\boldsymbol{\alpha}}$, with
  $\boldsymbol{\alpha}=\left(\alpha_1,\alpha_2\right)$,
  $\abs{\boldsymbol{\alpha}}=\alpha_1+\alpha_2\leq k-2$, such that
  \begin{equation*}
    \forall (x,y) \in E,\quad m_{\boldsymbol{\alpha}}(x,y) \defeq
    \frac{(x-x_E)^{\alpha_1}(y-y_E)^{\alpha_2}}
    {h_E^{\alpha_1+\alpha_2}}\,,
  \end{equation*}
  being $(x_E,y_E)$ the centroid of the cell and $h_E$ its diameter.
\end{enumerate}
For any element $E\in\Th$, given a function $v_\delta\in\Vh^E$, it can be seen
\cite{Beirao2013a,Beirao2014d} that the values of its degrees of
freedom are uniquely defined by its $\lebl{E}$-orthogonal projection on
$\Poly{k-1}{E}$, denoted by $\proj{k-1,E}{}v_\delta$, and the
orthogonal projection of its gradient on
$\Poly{k-1}{E} \times \Poly{k-1}{E}$, denoted by
$\proj{k-1,E}{}\nabla v_\delta$. A basis of the local VEM space is defined
implicitly as the set of functions that are Lagrangian with respect to
the degrees of freedom.

To discretize \eqref{eq:darcy-model} by the Virtual Element method we
suppose to know, for each $E\in\Th[i]$, $i= 1,\ldots,N$, a bilinear
form $\vemstab[]{}{}\colon \Vh^E \times \Vh^E \to \mathbb{R}$ such
that, $\forall v_\delta \in \Vh^E \cap \ker \proj[\nabla]{k,E}{}$,
\begin{equation*}
  \label{eq:SEscales}
  \exists c_\ast\,,c^\ast >0 \colon
  c_\ast\scal[E]{\K_i \nabla v_\delta}{\nabla v_\delta} \leq \vemstab[E]{v_\delta}{v_\delta}
  \leq c^\ast \scal[E]{\K_i \nabla v_\delta}{\nabla v_\delta} \,.
\end{equation*}
With the above assumption, we define, for each $E\in\Th[i]$,
$i = 1,\ldots,N$, the bilinear form
$\ahE{}{}\colon \Vh^E \times \Vh^E \to \mathbb{R}$ such that,
$\forall u_\delta,v_\delta\in \Vh$,
\begin{equation*}
  \ahE{u_\delta}{v_\delta} = \scal[E]{\K_i\Pi_{k-1}^0\nabla u_\delta}{\Pi^0_{k-1}
    \nabla v_\delta} + \K_i
  \vemstab[E]{\left(I-\Pi^\nabla_k\right)u_\delta}{\left(I-\Pi^\nabla_k\right)v_\delta}
  \,,
\end{equation*}
and the global bilinear form
$\ah{}{}\colon \Vh^D \times \Vh \to \mathbb{R}$ such that
\begin{equation*}
  \ah{u_\delta}{v_\delta} = \sum_{i = 1}^N \sum_{E\in\Th[i]}
  \ahE{\eval{u_\delta}{F_i}}{\eval{v_\delta}{F_i}}.
\end{equation*}

Finally, the Virtual Element discretization of \eqref{eq:darcy-model}
is: find $H_\delta \in \Vh^D$ such that, $\forall v_\delta \in \Vh$
\begin{multline}
  \label{eq:darcy-discrmodel}
  \ah{H_\delta}{v_\delta} = \sum_{i=1}^N \left( \sum_{E\in\Th[i]}
    \scal[E]{f_i }{\proj{k-1,E}{\eval{v}{F_i}}} \right)
  \\
  + \dual[\sobh{-\frac12}{\Gamma_{N_i}}, \sobh{\frac12}{\Gamma_{N_i}}]
  {g_{i}^{N}}{\eval{v}{\Gamma_{N_i}}} \,.
\end{multline}
We remark that the continuity conditions are automatically satisfied
by the definition of the functional space and the degrees of freedom,
viable because $\Th$ is a globally conforming discretization.

In \cite{BBapost}, a residual a posteriori estimate was derived for
the Laplace problem proving the equivalence between the estimator
and the error with respect to a suitable polynomial projection of the VEM solution. The extension of this estimate to the
case of a globally conforming discretization of the Laplace problem on
a DFN is quite straightforward. Let us define the following error
measure:
\begin{equation*}
  \ennorm{v} = \sup_{w\in V}\frac{\sum_{i=1}^N
    \a[i]{v_i}{w}}{\left(\sum_{i=1}^N\norm[F_i]
      {\sqrt{\K}\nabla w}^2\right)^{\frac12}} \,,
\end{equation*}
then, we denote $H_\delta^\pi = \proj[\nabla]{k}{H_\delta}$ and define
\begin{multline*}
  est_\delta^2 = \sum_{i=1}^N\left( \sum_{E\in\Th[i]}
    \frac{h^2_E}{\K_i} \norm[E]{\proj{k-1}{}f+\K_i\Delta H_{\delta
        i}^\pi }^2 + \sum_{e\in\Ehint[i]} \frac{h_e}{\K_i}
    \norm[e]{\jmpnormder{H_{\delta i}^\pi}{e}{\K_i}}^2 + \right.
  \\
  \left. +\sum_{e\in\EhN[i]} \frac{h_e}{\K_i}
    \norm[e]{\jmpnormder{H_{\delta i}^\pi}{e}{\K_i} - g_{i}^N}^2 +
    \sum_{E\in\Th}\frac{h^2_E}{\K_i}\norm[E]{f-\proj{k-1}{}f}^2
  \right)
  \\
  + \sum_{\substack{S\in\CAL{S}\colon\\ \I(S) = \{i,j\}}} \left(
    \sum_{e\in\Eh[S]} \frac{h_e}{\min \{\K_i,\K_j\}}
    \norm[e]{\jmpnormder{H_{\delta i}^\pi}{e}{\K_i} +
      \jmpnormder{H_{\delta j}^\pi}{e}{\K_j}}^2 \right)\,,
\end{multline*}
where $\forall i=1,\ldots N$,
$\Ehint[i]$ is the set of edges of $\Th[i]$ such that,
$\forall e \in\Ehint[i]$, $e\cap S = \emptyset$ $\forall S\in\CAL{S}_i$ and
$e\cap \Gamma_{N_i} = \emptyset$, $\EhN[i]$ is the set of edges $e$ of
$\Th[i]$ such that $e \cap \Gamma^N_i \neq \emptyset$ and
$\forall S\in\CAL{S}$, $\Eh[S]$ is the set of edges $e$ of $\Th$ such
that $e \cap S \neq \emptyset$. Then, there exist two constants
$c,C>0$ independent of $\delta$ such that
\begin{equation}
  \label{eq:estim-equivalence}
  c \cdot est_\delta \leq \ennorm{H-H_\delta^{\pi}} \leq C \cdot est_\delta \,.
\end{equation}
In view of an adaptive approach, we define, for each cell $E$, a local
estimator $est_{\delta,E}$, such that the estimators defined on the edges are
split among neighbouring cells according to their areas.


\section{DFN Minimal mesh construction}
\label{sec:minimal-mesh}
In this section we introduce the strategy used for the construction of the initial coarse polygonal mesh on the DFNs. This mesh is obtained by the construction of convex polygons representing sub-fractures, i.e. portion of fractures not crossed by traces that can have traces or portion of traces or extension of traces only on the boundary. Given a fracture and the set of its traces there exist many partitions of the fracture in sub-fractures with a different number of sub-fractures and different quality of the element produced. The construction of the minimal mesh is a complex trade off between the number of the elements and the quality of the elements produced. In our approach we aim at limiting the number of elements. Improvement of the mesh quality is transferred to a following suitable refining strategy.

The approach we follow is an iterative splitting of the leaves of a tree structure. We start from the original fracture that is the root cell of the structure, then we select a trace and we split the given cell along the trace or an extension of the trace producing two children cells.
Then we proceed iteratively choosing a new trace and cutting the leaves of the tree with the selected trace, exclusively in the case the trace is intersecting the internal part of the cell. At each iteration each cell can be split in two children cells or can be modified in a unique child cell that is the same polygon with one of the edges split in two aligned sub-edges in order to guarantee conformity of the global mesh. The algorithm is sketched in Algorithm~\ref {alg:MinimalMesh}.

In order to control the number of cells produced by the algorithm we cut the current cells with a suitable order of traces. We start considering the traces that cross the fracture intersecting two boundary edges of the fracture. Then we continue considering the remaining traces from the longest to the shortest. Considering the traces in this order usually yields to a smaller number of cells.

\begin{algorithm}
    \caption{Minimal mesh}
    Given a fracture and the set of traces create a tree structure with the fracture as root cell
    \begin{algorithmic}[1]
      \FOR{All the traces}
      \FOR{All the leaves cells of the tree}
      \STATE{Compute the intersection of the trace with the cell}
      \IF{There is an intersection with the internal part of the cell}
       \STATE Split the cell in two children cells
       \ENDIF
       \STATE Update the neighbouring cell with the new edges
       \ENDFOR
       \ENDFOR
    \end{algorithmic}
    \label{alg:MinimalMesh}
  \end{algorithm}

\section{Refinement and Marking algorithms}
\label{sec:algorithm}

In this section we briefly introduce the algorithms used for marking the cells with largest estimators to be refined in order to reduce the discretization error and the different algorithms tested for the refinement of polygonals cells. This refinement step is part of the usual refinement process SOLVE-ESTIMATE-MARK-REFINE \cite{Doerfler}.

\subsection{Marking Strategy}
The marking strategy (see Algorithm~\ref{alg:Marking}) of the cells to be refined is simply based on the selection of all the cells with the largest error estimators.
We mark the cells starting from those with largest estimators $est_{\delta,E}^2$ up to when the cumulative error estimator $est_{\delta,\Th[{mark}]}^2$ of the marked cells is a given ratio $C$ of the total error estimator $est_{\delta,\Th}^2=\sum_{E\in\Th} est_{\delta,E}^2$. In this algorithm we accept the sorting cost for the estimators vector \cite{Doerfler} in order to maximally contain the refinement iterations that, in practical applications to large scale DFN simulations, can be quite expensive in the last refinement steps.

  \begin{algorithm}
    \caption{Cells marking algorithm}
    Given a convex polygon
    \begin{algorithmic}[1]
      \STATE Compute the cell error estimators $est_{\delta,E}^2$ and save them in a vector $Vest_{\delta}^2$.
      \STATE Compute the total error estimator $est_{\delta,\Th}^2$
      \STATE Sort the cell estimators vector $Vest_{\delta}^2$
      \STATE Choose a constant $0 < C < 1$
      \STATE{$i=0$}
      \STATE{$est_{\delta,\Th[{mark}]}^2=est_{\delta,\Th[{mark}]}^2+Vest_{\delta}^2(i)$}
      \REPEAT
      \STATE{$i++$}
      \STATE{$est_{\delta,\Th[{mark}]}^2=est_{\delta,\Th[{mark}]}^2+Vest_{\delta}^2(i)$}
      \UNTIL{$est_{\delta,\Th[{mark}]}^2 \leq  C * est_{\delta,\Th}^2$}
      \STATE Mark the cells corresponding to the first $i$ positions of $Vest_{\delta}^2$
    \end{algorithmic}
    \label{alg:Marking}
  \end{algorithm}

\subsection{Refinement algorithms}
In this section we introduce four different refinement algorithms used for cutting marked convex cells in two convex sub-cells. All the algorithms are based on a similar approach and differ for the choice of the cutting direction. The common approach is described in the Algorithm~\ref{alg:Refinement} and the four different approaches differ for the Step~\ref{alg:Direction} of Algorithm~\ref{alg:Refinement}.

\begin{algorithm}
    \caption{Refinement algorithm for convex polygons}
    Given a convex polygon
    \begin{algorithmic}[1]
      \STATE Check aspect ratio (AR) \label{alg:Check}
      \STATE Compute the centroid $\mathbf{X}_c$
      \STATE Choose the cutting direction \label{alg:Direction}
      \STATE Build a straight line passing from the centroid with the chosen direction
      \STATE Choose the collapsing tolerance CollapseToll
       \FOR{Each edges of the cell}
       \IF{There is an intersection}
        \IF{The intersection point is near to the begin/end point of the edge according to the chosen tolerance CollapseToll} \label{alg:Direction-Modification}
        \STATE Change the cutting direction and set the intersection point as the begin/end point of the edge
        \ELSE
        \STATE The intersection point will be a new point in the mesh
        \ENDIF
       \ENDIF
       \ENDFOR
       \STATE Create the two children cells
       \STATE Update the neighbourhood
    \end{algorithmic}
    \label{alg:Refinement}
  \end{algorithm}


\begin{figure}
  \begin{subfigure}{0.46\linewidth}
    \includegraphics[width=\linewidth]{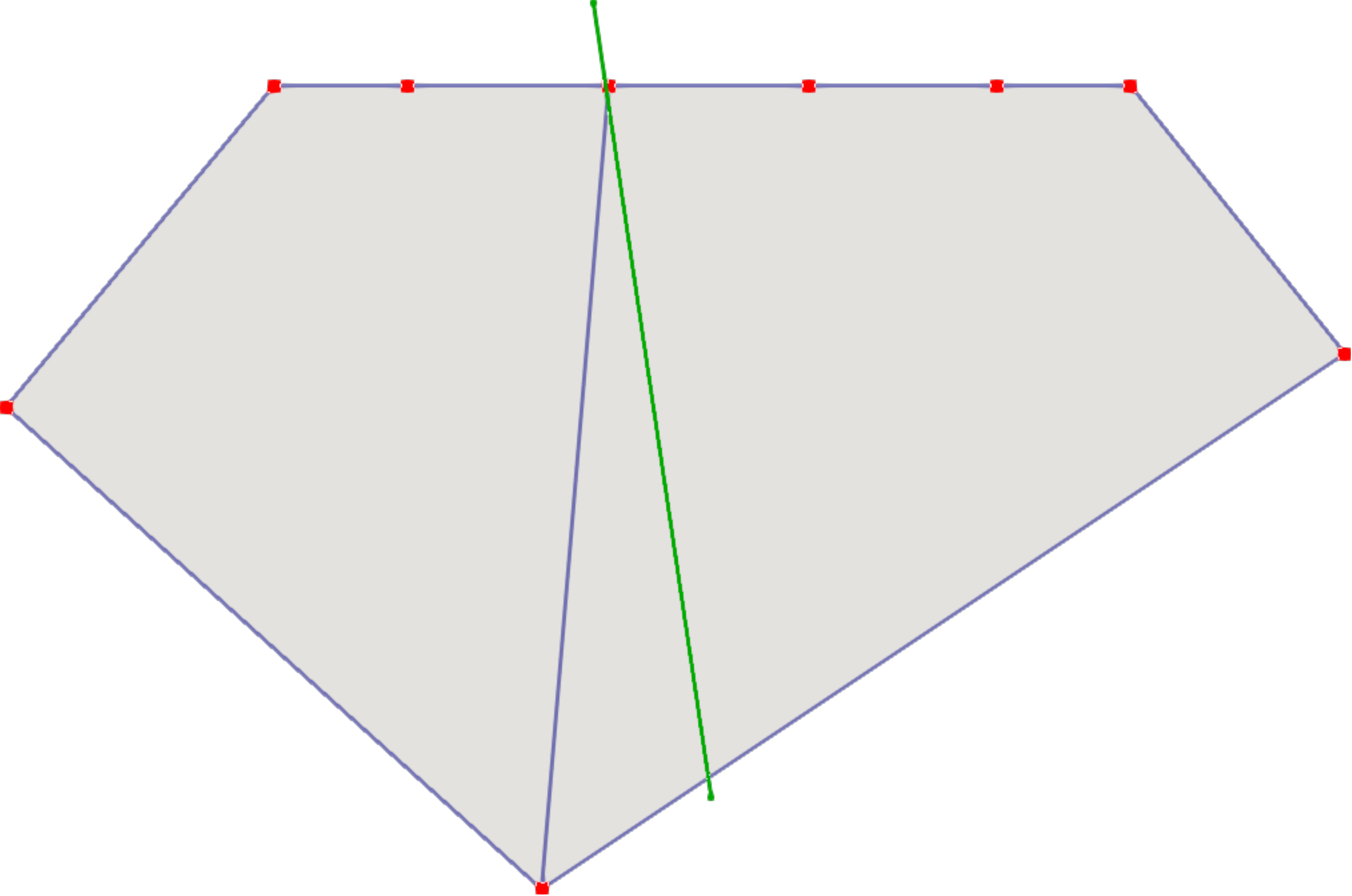}
    \caption{Maximum Momentum}
    \label{fig:Direction-MaxMom}
  \end{subfigure}
  \hfil
  \begin{subfigure}{0.46\linewidth}
    \includegraphics[width=\linewidth]{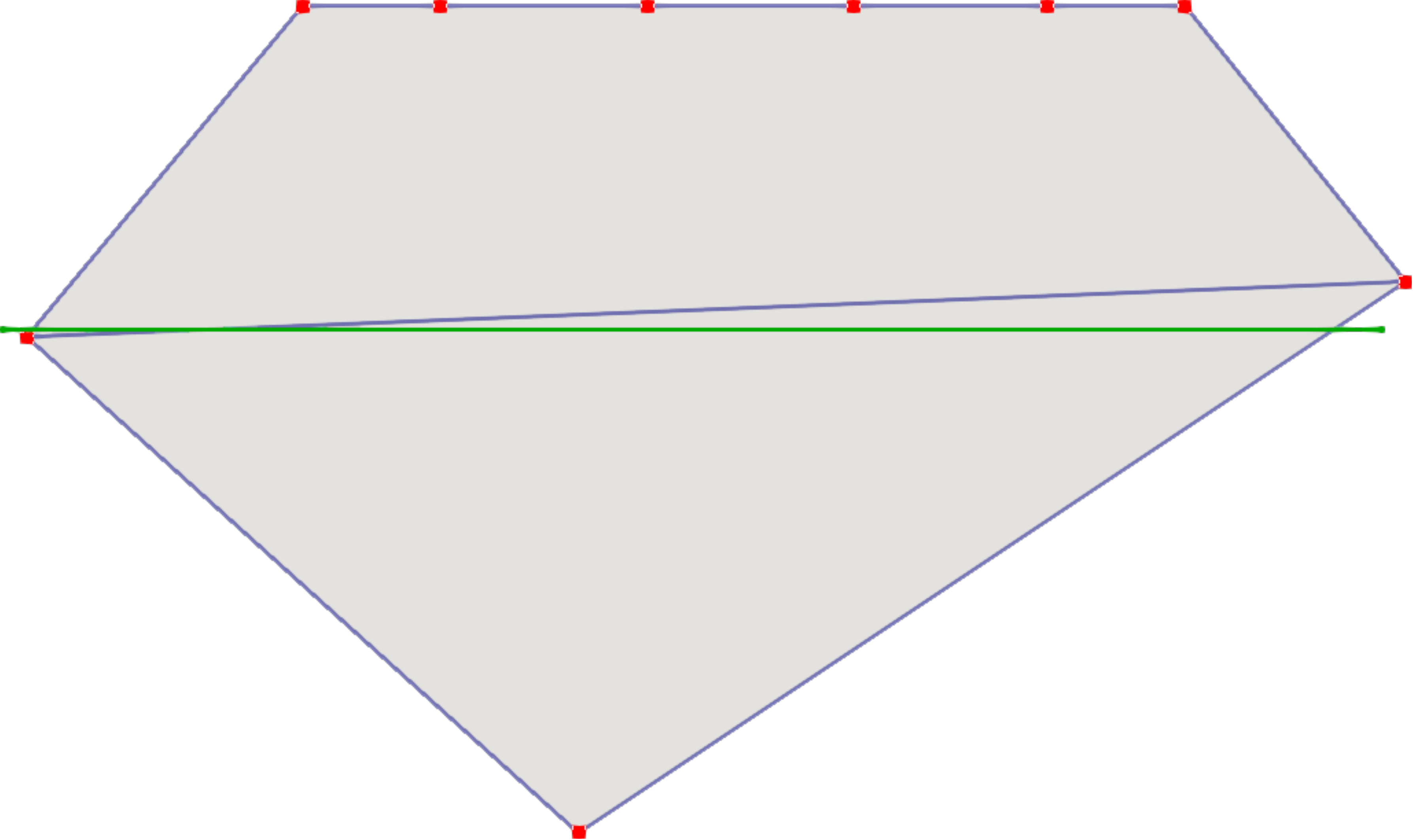}
    \caption{Trace Direction}
    \label{fig:Direction-TrDir}
  \end{subfigure}
  \begin{subfigure}{0.46\linewidth}
    \includegraphics[width=\linewidth]{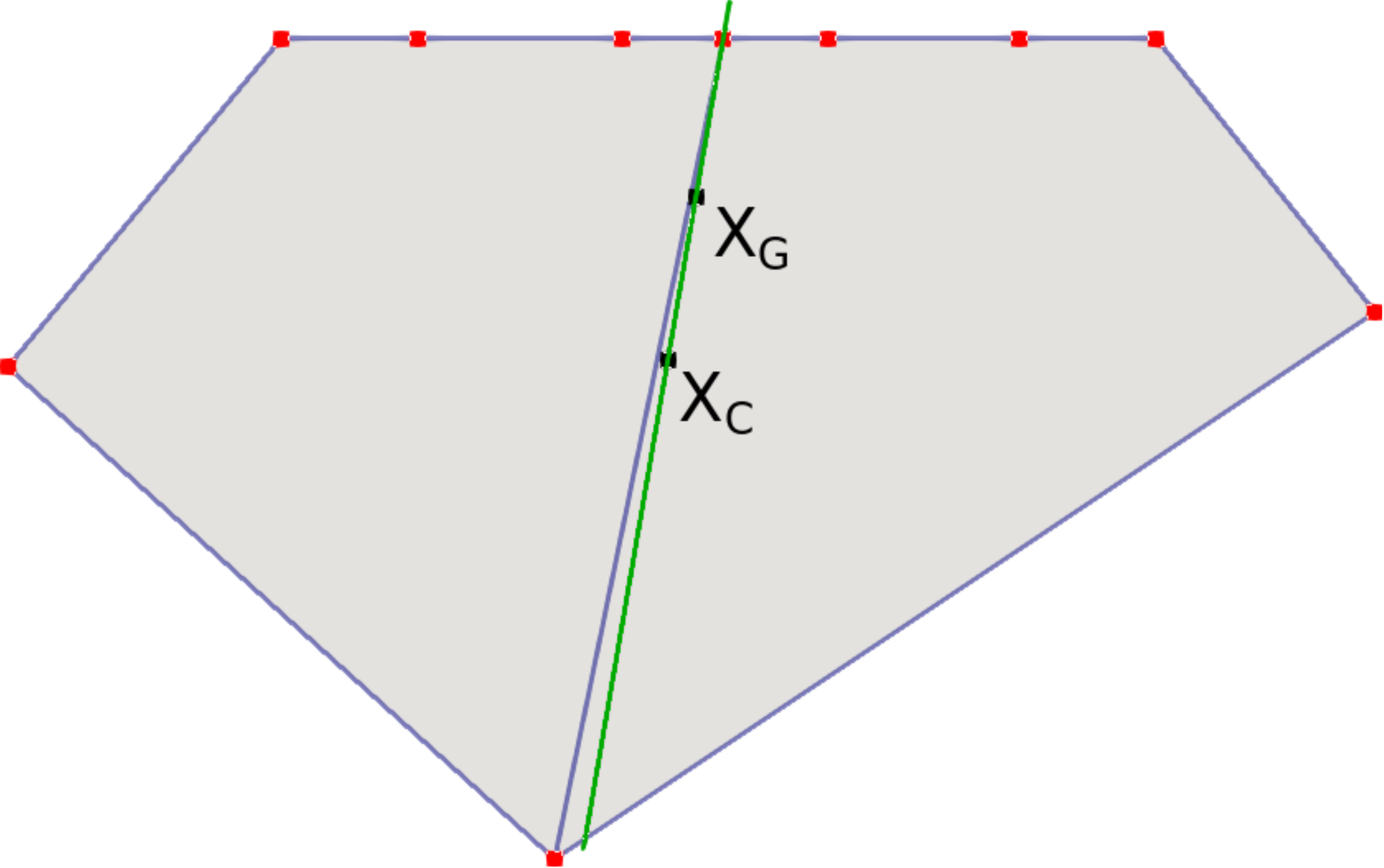}
    \caption{Maximum Number of Points}
    \label{fig:Direction-MaxPnt}
  \end{subfigure}
  \hfil
  \begin{subfigure}{0.46\linewidth}
    \includegraphics[width=\linewidth]{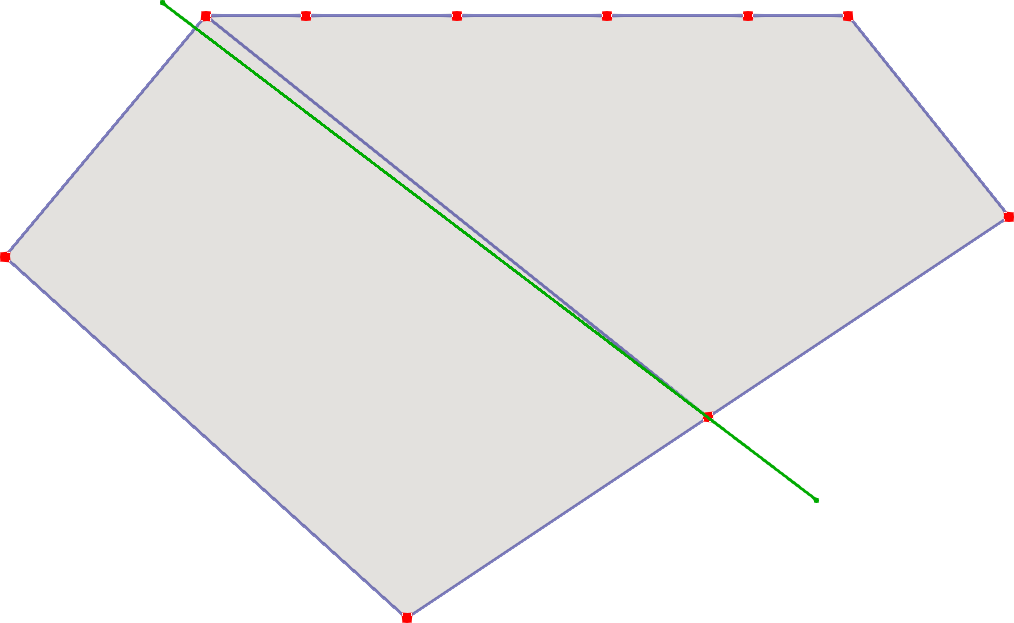}
    \caption{Maximum Edge}
    \label{fig:Direction-MaxEdg}
  \end{subfigure}
  \caption{Different strategies for the selection of the cutting direction.}
  \label{fig:Direction}
\end{figure}


In the following we compare four different refinement options for choosing the cutting direction denoted Maximun Momentum (MaxMom), Trace Direction (TrDir), Maximum Number of Points (MaxPnt) and Maximum Edge (MaxEdg). At the Step~\ref{alg:Direction-Modification} of the Algorithm~\ref{alg:Refinement} we accept to slightly modify the chosen cutting direction in order to avoid the proliferation of edges and vertices of the new cells and consequently of VEM degrees of freedom that do not efficiently increase the quality of the solution and to avoid the generation of very small edges that could sometimes induce stability problems \cite{Beirao2017,BBorth}.

We collapse the points given by the intersection of the cutting direction with an edge to existing vertices of the cell when the distance of the intersection from the closest vertex are smaller than the tolerance CollapseToll multiplied by the length of the edge intersected.

For some of the proposed cutting directions we switch from the selected refinement criterion to the Maximum Momentum criterion in order to avoid the generation of cells with a huge aspect ratio (AR) defined as the ratio between the longest distance between the centroid and the vertices and the smallest distance between the centroid and the edges. In order to avoid the generation of elements with a huge aspect ratio, we define a fixed value denoted by MaxAR and if the aspect ratio of the cell to cut is larger than MaxAR, the cutting criterion is always MaxMom.

\subsubsection{Maximum Momentum (MaxMom)}
This cutting strategy of marked cells is based on the choice of cutting direction that is orthogonal to the direction of eigenvector associated to the smallest eigenvalue of the inertia tensor of the cell (Figure~\ref{fig:Direction-MaxMom}).

\subsubsection{Trace Direction (TrDir)}
This cutting strategy is based on the choice as cutting direction of the direction parallel to a trace (if any), see Figure~\ref{fig:Direction-MaxPnt}, where we have assumed that the edge with several aligned points belongs to the trace. The rationale of this approach is related to the known property of the solution to the considered problem that displays strong gradient components in the direction orthogonal to the traces. If a cell intersects more than one trace we switch to MaxMom criterion.

\subsubsection{Maximum Number of Points (MaxPnt)}
This cutting strategy is based on the choice as cutting direction of the direction of the vector connecting the centroid of the cell $X_C$ to the center of mass of the vertices $X_G$, see Figure~\ref{fig:Direction-MaxPnt}. This vector points towards the region of the marked cell with the highest density of vertices and should split the cell balancing the vertices of the two new subcells. In this strategy it is mandatory to define the option MaxNP, to set how many points the cells can have.
This strategy switches to the MaxMom strategy in two cases: the first when the number of points of the cell to cut is less than MaxNP, the second when the centroid and the center of mass have a distance under a fixed tolerance.
We remark that this refinement strategy can be considered as a simple improvement of the MaxMom strategy being the refinement strategy different only for the cells with a large number of cells and considering that when the number of vertices of the cell is larger than MaxNP this refinement strategy aims at dropping the number of vertices of the two produced cells under MaxNP.

\subsubsection{Maximum Edge (MaxEdg)}
This cutting strategy is based on the choice as cutting direction of the direction that cuts the longest edge of the cell in half, see Figure~\ref{fig:Direction-MaxEdg}, where the longest edge is the right-bottom edge. If a cell displays aligned edges, these are considered as one unique edge.


\section{Numerical Results: optimality and effectivity index}
\label{sec:num_res_TP}

In order to validate our refinement algorithm we test it on two simple
DFNs for which an exact solution is known. We consider two DFNs with
two and three fractures, labelled as Problem 1 and Problem 2,
respectively.

This first set of tests aims at validating the equivalence relation
stated in \eqref{eq:estim-equivalence} between the error and the error
estimator. We have tested this equivalence relation on the meshes
produced during the adaptive process although this property holds true
on any sufficiently refined mesh.  We apply our adaptive algorithm and
compare at each refinement iteration the error and the error estimator
computing the {\sl effectivity index}
\begin{equation*}
  \varepsilon = \frac{err}{est_\delta}
\end{equation*}
in order to verify that it is independent of the mesh size obtained by
adaptive refinements. See \cite{BBapost} for the same analysis
performed on uniformly refined meshes.

As stopping criterion for the adaptive process, we require the
following condition on the estimated relative error:

\begin{equation}
  \label{eq:relerr_stopping}
  \frac{est_\delta}{\left(\sum_{i=1}^N
      \norm[F_i]{\sqrt{\K_i}\nabla \proj[\nabla]{k}{}H_{\delta i}}^2
    \right)^{\frac12}} \leq 0.05 \,.
\end{equation}

All the simulations here presented are performed with the following
methods and parameters: VEM orders from $1$ to $4$, Preconditioned Conjugate
Gradient as linear solver with relative stopping residual $1.0e-15$,
the preconditioner being an incomplete Cholesky factorization implemented
as described in \cite{EigPCG}, CollapseToll = 0.2, MaxAR = 10,
MaxNP = 12 and C = 0.50.

\subsection{Problem $1$}

\begin{figure}
  \centering
  \begin{subfigure}{0.46\linewidth}
    \includegraphics[width=\linewidth]{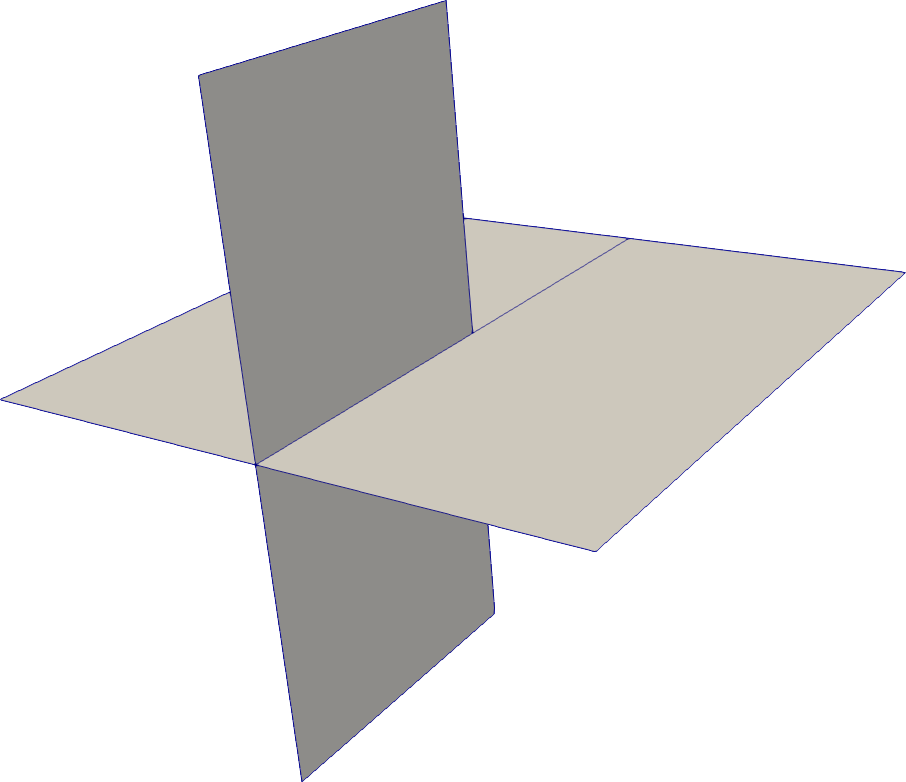}
    \caption{Initial mesh}
    \label{TipEnrich:DFN-mesh-1}
  \end{subfigure}
  \hfill
  \begin{subfigure}{0.46\linewidth}
    \includegraphics[width=\linewidth]{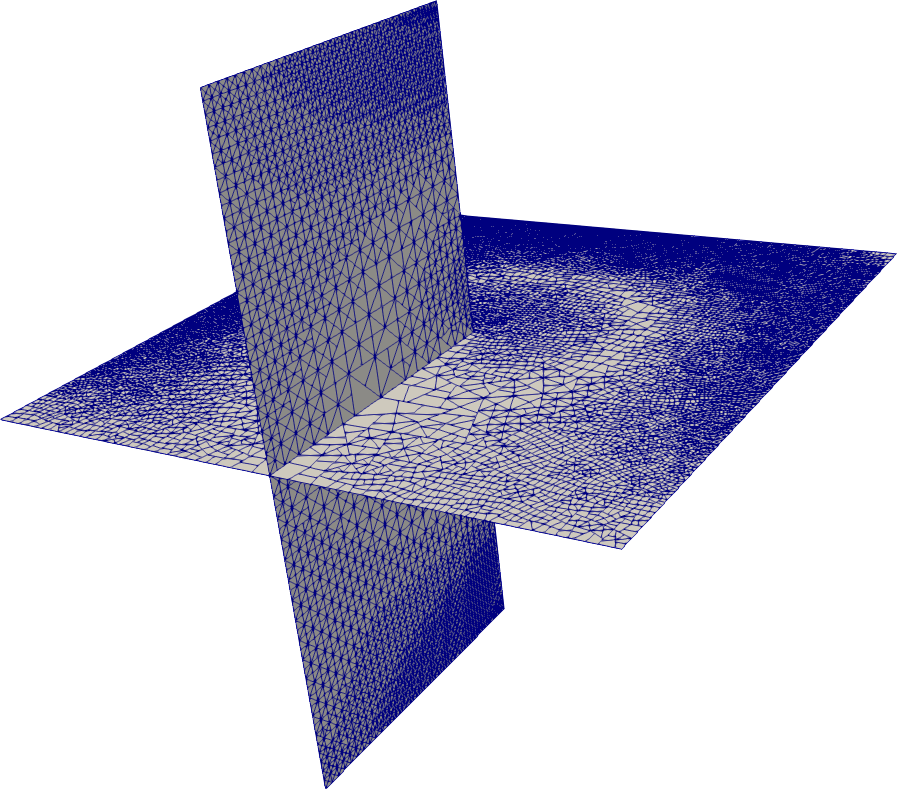}
    \caption{Step 30 MaxMom}
    \label{TipEnrich:DFN-mesh-30MM}
  \end{subfigure}
  \hfill
  \begin{subfigure}{0.46\linewidth}
    \includegraphics[width=\linewidth]{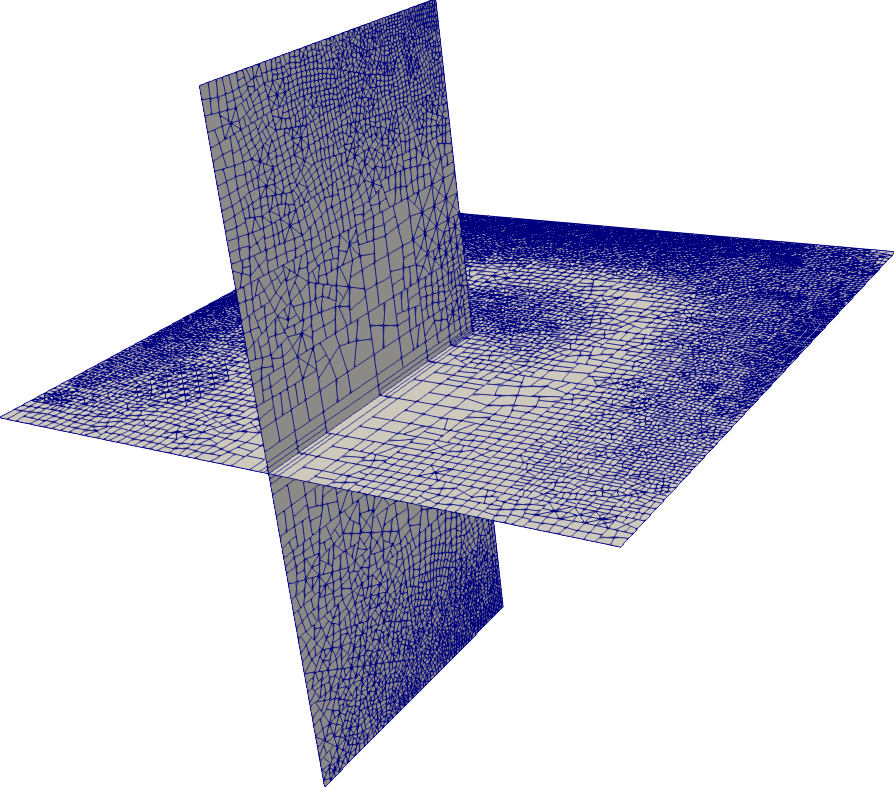}
    \caption{Step 30 TrDir}
    \label{TipEnrich:DFN-mesh-30TD}
  \end{subfigure}
  \hfill
  \begin{subfigure}{0.46\linewidth}
    \includegraphics[width=\linewidth]{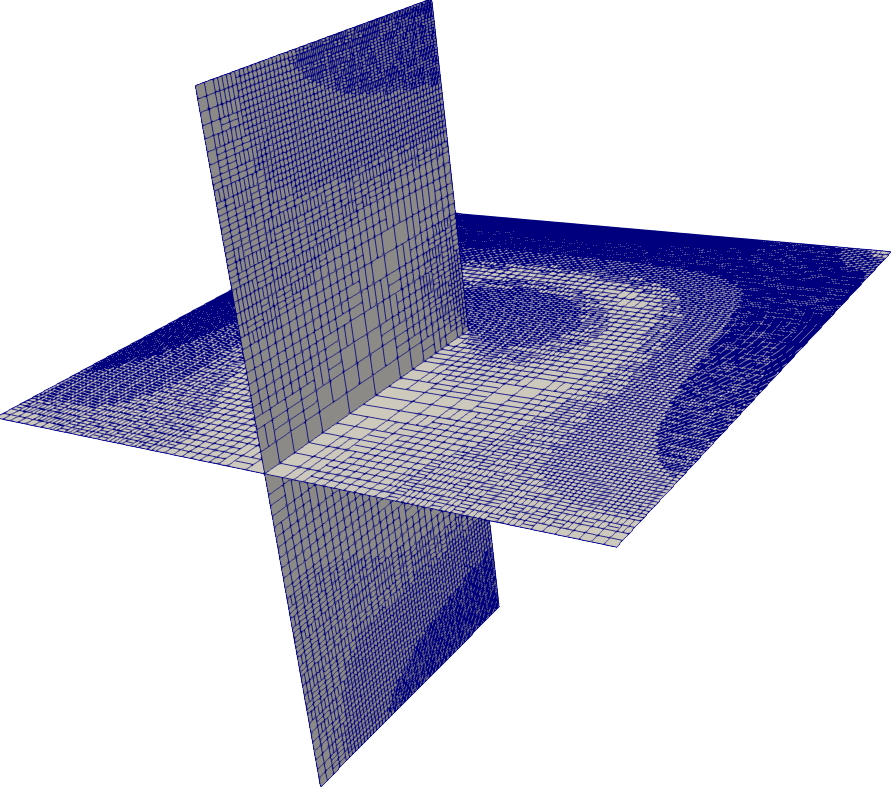}
    \caption{Step 30 MaxEdg}
    \label{TipEnrich:DFN-mesh-30ME}
  \end{subfigure}
  \caption{Problem 1: DFN with meshes at refining step 30.}
  \label{fig:TipEnrich-mesh}
\end{figure}
\begin{figure}
  \centering
  \begin{subfigure}[b]{0.49\linewidth}
    \includegraphics[width=\linewidth]{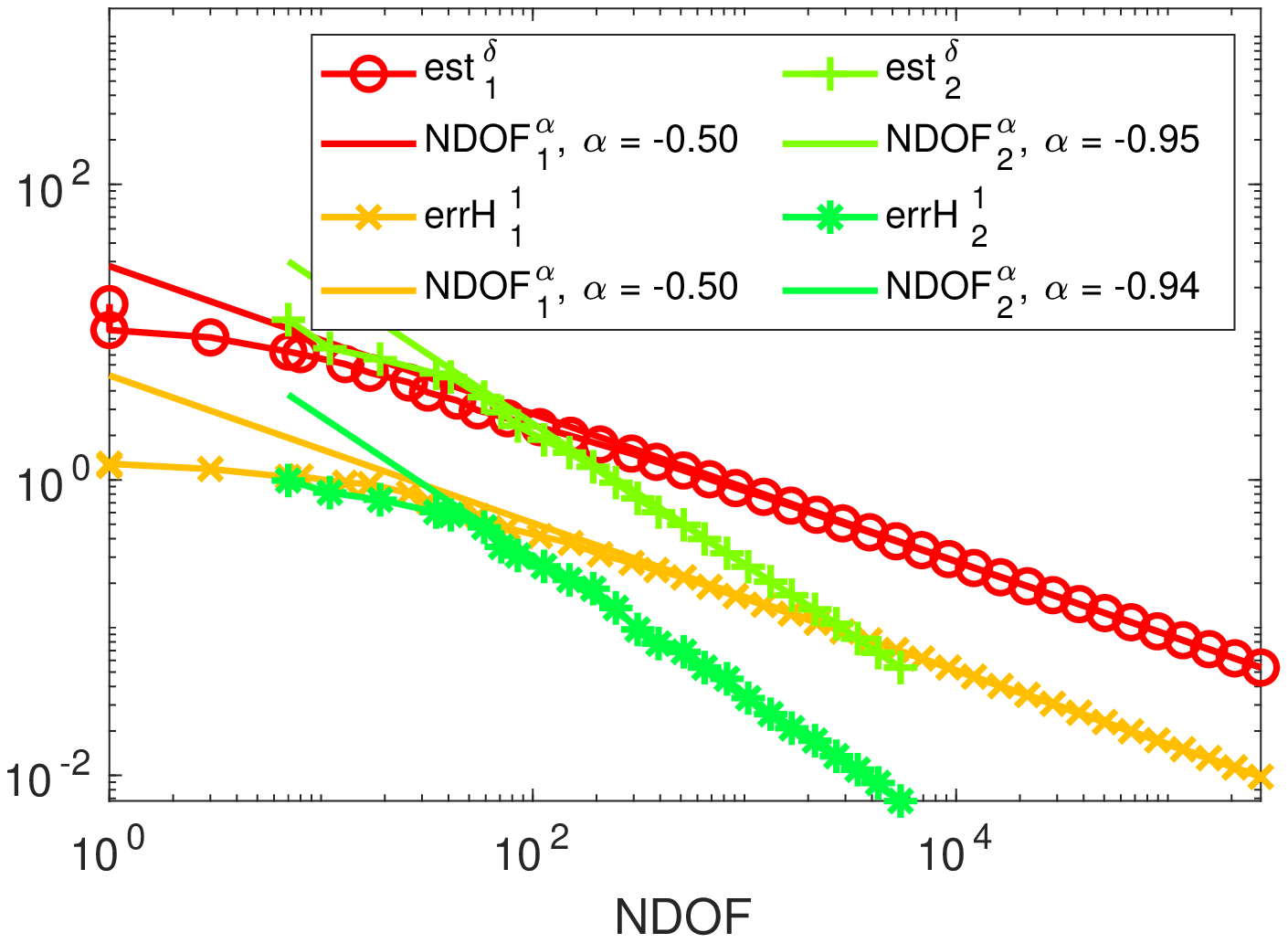}
    \label{fig:p1_MM}
  \end{subfigure}
  \hfill
  \begin{subfigure}[b]{0.49\linewidth}
    \includegraphics[width=\linewidth]{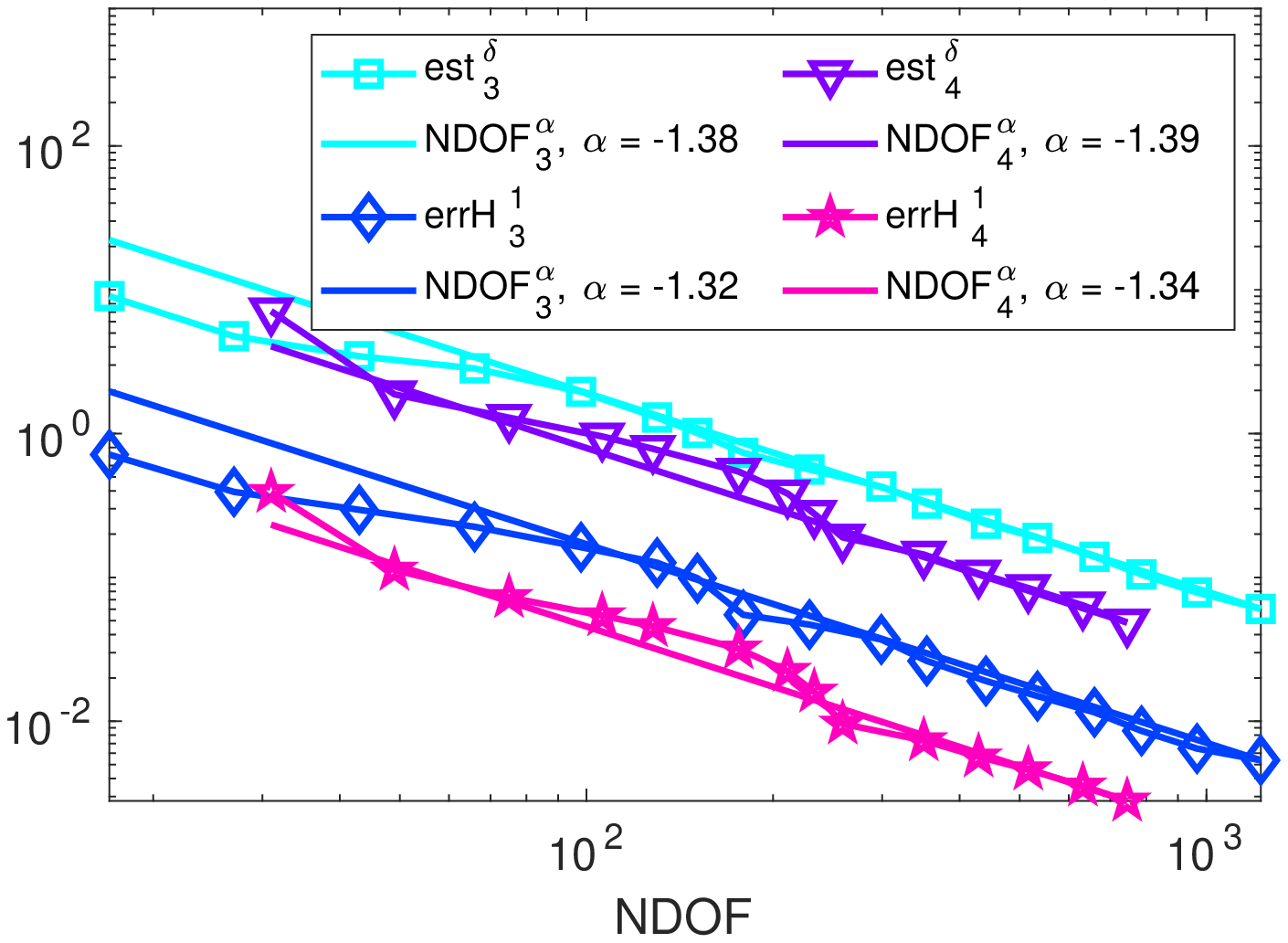}
    \label{fig:p1_TD}
  \end{subfigure}
  \caption{Problem 1: Rates of convergence ($\alpha$) for error and
    estimator with MaxMom refinement.}
  \label{fig:TipEnrich-Rates}
\end{figure}
\begin{table}[h]
  \centering \begin{tabular}{c|c|c|c|c}
TypeRef & Order 1 & Order 2 & Order 3 & Order 4 \\ 
\hline 
MaxMom & $\alpha = -0.50$ & $\alpha = -0.95$ & $\alpha = -1.38$ & $\alpha = -1.39$ \\ 
   & $\alpha_{err} = -0.50$ & $\alpha_{err} = -0.94$ & $\alpha_{err} = -1.32$ & $\alpha_{err} = -1.34$ \\ 
TrDir & $\alpha = -0.49$ & $\alpha = -0.93$ & $\alpha = -1.67$ & $\alpha = -1.56$ \\ 
   & $\alpha_{err} = -0.50$ & $\alpha_{err} = -0.92$ & $\alpha_{err} = -1.63$ & $\alpha_{err} = -1.30$ \\ 
MaxEdg & $\alpha = -0.50$ & $\alpha = -0.93$ & $\alpha = -1.37$ & $\alpha = -1.55$ \\ 
   & $\alpha_{err} = -0.50$ & $\alpha_{err} = -0.93$ & $\alpha_{err} = -1.41$ & $\alpha_{err} = -1.38$ \\ 
\end{tabular}
  \caption{Problem 1: Rates of convergence for the estimator ($\alpha$) and the error ($\alpha_{err}$) with refinement criteria MaxMom, TrDir and MaxEdg.}
  \label{tab:TipEnrich-Rates2}
\end{table}
\begin{table}
  \centering
  \begin{tabular}{c|c|c|c|c}
Step & $\#$ Cells & NDOF & $\varepsilon$ & PCG-It \\ 
\hline 
1 & 4 & 1 & 0.0064 & 1 \\ 
5 & 16 & 8 & 0.0221 & 5 \\ 
10 & 58 & 44 & 0.0318 & 12 \\ 
15 & 232 & 208 & 0.0326 & 23 \\ 
20 & 958 & 910 & 0.0345 & 43 \\ 
25 & 3997 & 3901 & 0.0343 & 88 \\ 
30 & 16696 & 16214 & 0.0337 & 181 \\ 
35 & 69827 & 67074 & 0.0334 & 370 \\ 
39 & 219297 & 207112 & 0.0333 & 662 \\ 
40 & 291581 & 274982 & 0.0332 & 777 \\ 
\end{tabular}
  \caption{Problem 1: MaxMom refinement Order 1.}
  \label{tab:TipEnrich-MM}
\end{table}
\begin{table}
  \centering
  \begin{tabular}{c|c|c|c|c}
Step & $\#$ Cells & NDOF & $\varepsilon$ & PCG-It \\ 
\hline 
1 & 4 & 1 & 0.0064 & 1 \\ 
5 & 14 & 5 & 0.0130 & 1 \\ 
10 & 46 & 34 & 0.0441 & 12 \\ 
15 & 167 & 147 & 0.0364 & 21 \\ 
20 & 672 & 712 & 0.0356 & 46 \\ 
25 & 2760 & 3232 & 0.0377 & 94 \\ 
30 & 11528 & 13848 & 0.0357 & 188 \\ 
35 & 48362 & 57738 & 0.0353 & 378 \\ 
39 & 152185 & 178922 & 0.0349 & 695 \\ 
40 & 202798 & 237727 & 0.0347 & 846 \\ 
\end{tabular}
  \caption{Problem 1: TrDir refinement Order 1.}
  \label{tab:TipEnrich-TD}
\end{table}
\begin{table}
  \centering
  \begin{tabular}{c|c|c|c|c}
Step & $\#$ Cells & NDOF & $\varepsilon$ & PCG-It \\ 
\hline 
1 & 4 & 1 & 0.0064 & 1 \\ 
5 & 16 & 7 & 0.0271 & 3 \\ 
10 & 58 & 43 & 0.0347 & 12 \\ 
15 & 239 & 223 & 0.0381 & 21 \\ 
20 & 974 & 1007 & 0.0382 & 47 \\ 
25 & 4101 & 4420 & 0.0382 & 100 \\ 
30 & 17254 & 18980 & 0.0380 & 205 \\ 
35 & 72401 & 80561 & 0.0386 & 414 \\ 
37 & 128831 & 139214 & 0.0385 & 545 \\ 
\end{tabular}
  \caption{Problem 1: MaxEdg refinement Order 1.}
  \label{tab:TipEnrich-ME}
\end{table}

The geometry and the parameters of this test problem are described in
detail in \cite[Section 5.1]{BPSb}.  The DFN is composed by two
fractures that are planar rectangles defined as
\begin{align*}
  F_1 = (-1,1) \times (-1,1) \times \{0\} \,;
  \quad
  F_2 = (-1,0) \times \{0\} \times (-1,1) \,.
\end{align*}
Thus, the only trace of the DFN is
$\Gamma_1 = (-1,0) \times \{0\} \times \{0\}$. The differential model
is set up by choosing the following exact solutions:
\begin{align*}
  h_1(x,y,z) &= \left(x^2-1\right) \left(y^2-1\right) \left(x^2+y^2\right)
               \cos \left( \frac12 \mathrm{arctan2}(x,y) \right) \,,
  \\
  h_2(x,y,z) &= -\left(z^2-1\right) \left(x^2-1\right) \left(x^2+z^2\right)
               \cos \left( \frac12 \mathrm{arctan2}(z,x) \right) \,,
\end{align*}
where $\mathrm{arctan2}(x,y)$ is the four-quadrant inverse tangent,
that is the $\mathrm{arctan}$ of $y/x$ in $[-\pi,\pi]$. The
transmissivity of both fractures is set to $1$.

In Figure \ref{fig:TipEnrich-mesh} we report some of the meshes
generated during the refinement process.  We remark that in this test
problem the presence of a known forcing function on the fractures
induces a refinement in the whole domain, this will not be the same
for the flow simulations in which the rock matrix surrounding the DFNs
is considered impervious (Section \ref{sec:num_res_DFN}).

In Figure \ref{fig:TipEnrich-Rates} we report the convergence
behaviour of the error (errH$^1$) and of the error estimator
(est$_\delta$) and the final rates of convergence $\alpha$ with
respect to the total number of degrees of freedom (NDOF):
$est_\delta \sim \left(NDOF\right)^\alpha$, computed on the basis of
the last five refinement iterations, considering the strategy
MaxMom criterion.
We can clearly appreciate a parallel
behaviour of error and error estimator as well as the almost optimal
asymptotic rate of convergence very close to $-0.5$ and $-1$ for the VEM orders 1 and 2, respectively. For higher VEM orders the sub-optimal rates of convergence are due to the bounded Besov regularity of the solution around the internal trace-tip.
In Table \ref{tab:TipEnrich-Rates2} we report the rates of convergence
for the estimator and for the error obtained by the refinement strategies MaxMom, TrDir and MaxEdg.
For this problem we do not report results for the MaxPnt criterion that is always the MaxMom criterion being the number of vertices of the cells always smaller than MaxNP.

In Tables \ref{tab:TipEnrich-MM}, \ref{tab:TipEnrich-TD}
and \ref{tab:TipEnrich-ME} we report
the most significant quantities to describe the refinement process for
the MaxMom, TrDir and MaxEdge strategies: NCell is the total number of cells on
the DFN, NDOF is the total numbers of degrees of freedom,
$\varepsilon$ is the effectivity index, PCG-It is the number of
conjugate gradient iterations performed.  For all the strategies,
we highlight the relatively small variations
of $\varepsilon$ with respect to the large variations of number of
cells and degrees of freedom, after the first iterations corresponding
to very small number of NDOF.  We also remark the weak growing of
PCG-It with respect to the growing of NDOF.

\subsection{Problem $2$}
\begin{figure}
  \centering
  \begin{subfigure}{0.46\linewidth}
    \centering
    \includegraphics[width=\linewidth]{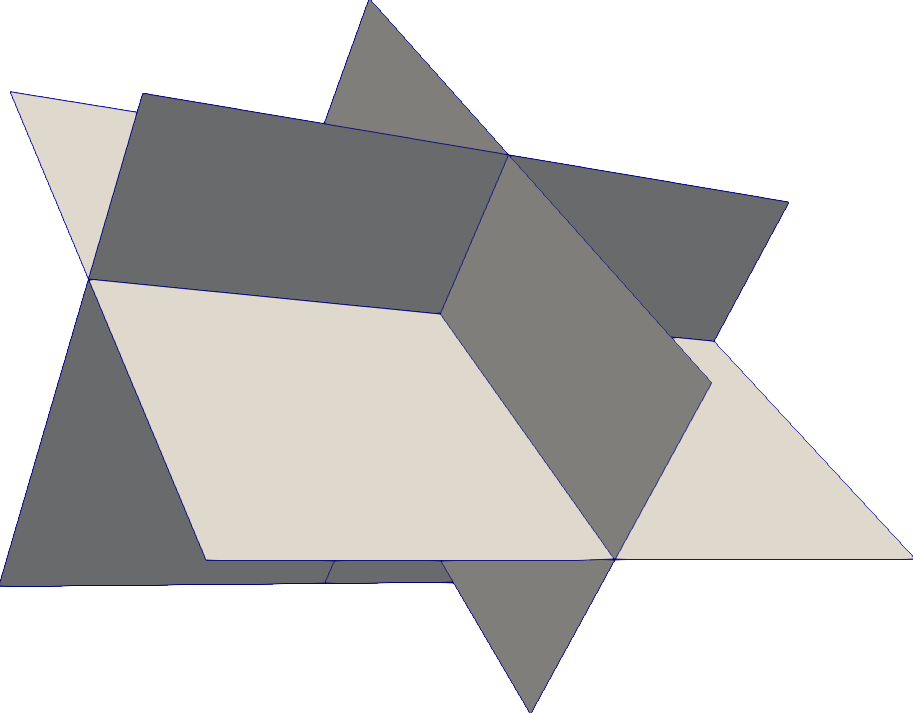}
    \caption{Initial mesh}
    \label{TipETut:DFN-mesh-1}
  \end{subfigure}
  \hfill
  \begin{subfigure}{0.46\linewidth}
    \includegraphics[width=\linewidth]{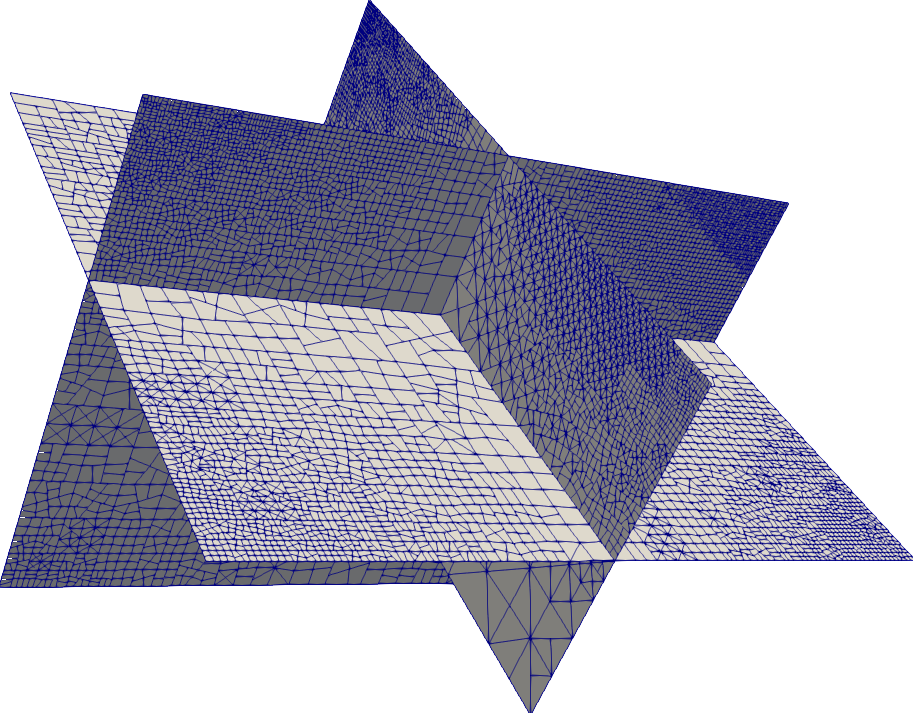}
    \caption{Step 25 MaxMom refinement Order 1}
    \label{TipETut:DFN-mesh-25MM}
  \end{subfigure}
  \hfill
  \begin{subfigure}{0.46\linewidth}
    \includegraphics[width=\linewidth]{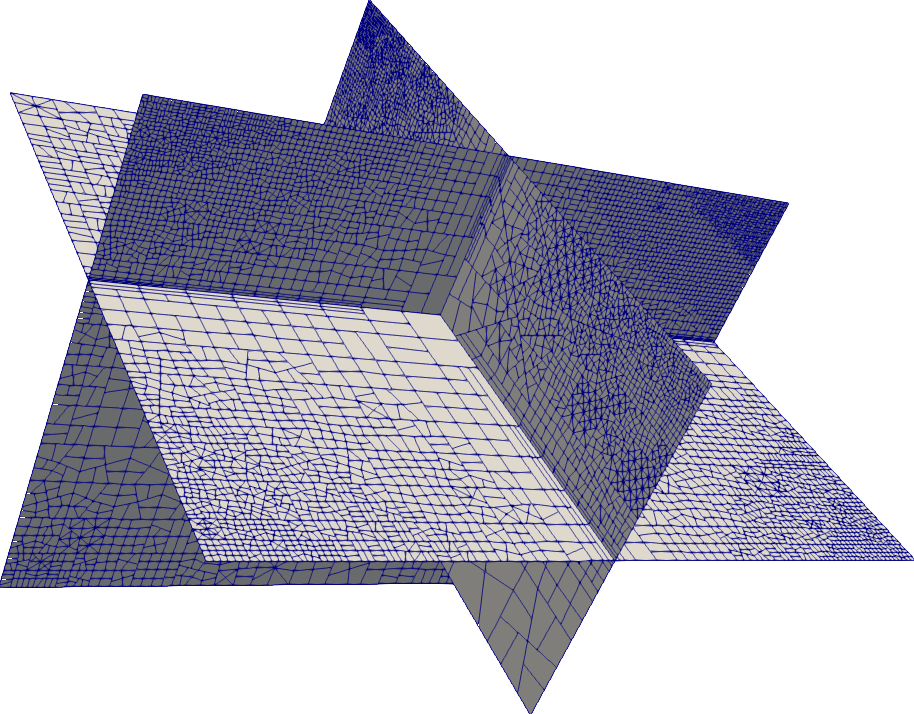}
    \caption{Step 25 TrDir refinement Order 1}
    \label{TipETut:DFN-mesh-25TD}
  \end{subfigure}
  \hfill
  \begin{subfigure}{0.46\linewidth}
    \includegraphics[width=\linewidth]{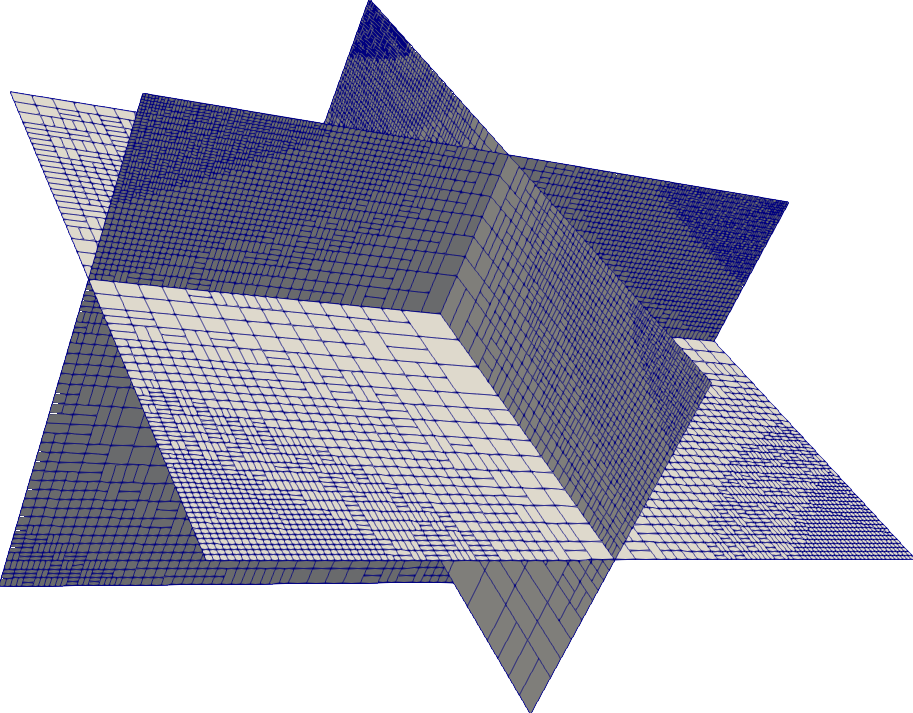}
    \caption{Step 25 MaxEdg refinement Order 1}
    \label{TipETut:DFN-mesh-25ME}
  \end{subfigure}
  \caption{Problem 2: DFN with meshes at refining step 25}
  \label{fig:TipETut_DFN-mesh}
\end{figure}
\begin{figure}
  \centering
  \begin{subfigure}[b]{0.49\linewidth}
    \includegraphics[width=\linewidth]{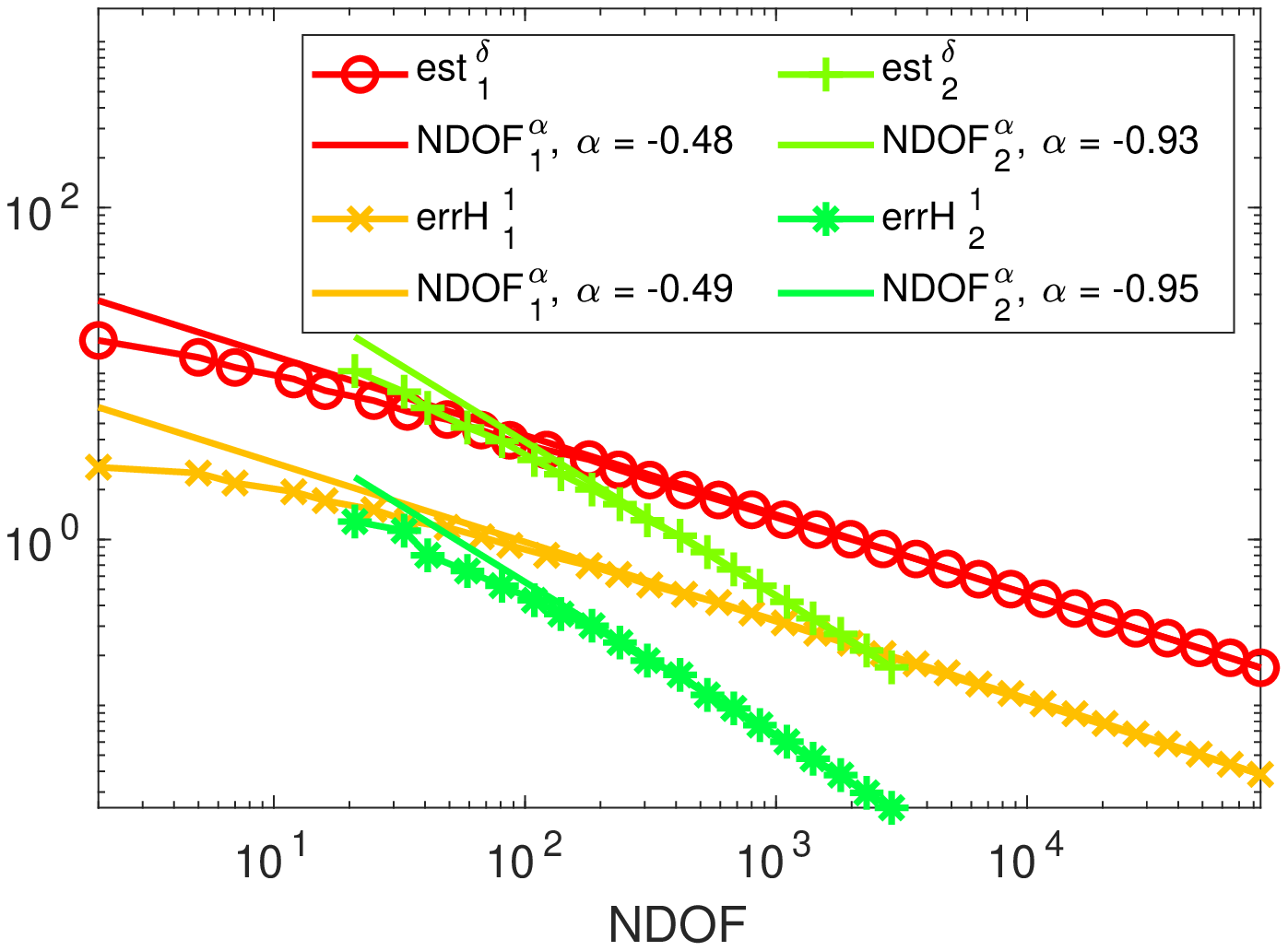}
    \label{fig:p2_MM}
  \end{subfigure}
  \hfill
  \begin{subfigure}[b]{0.49\linewidth}
    \includegraphics[width=\linewidth]{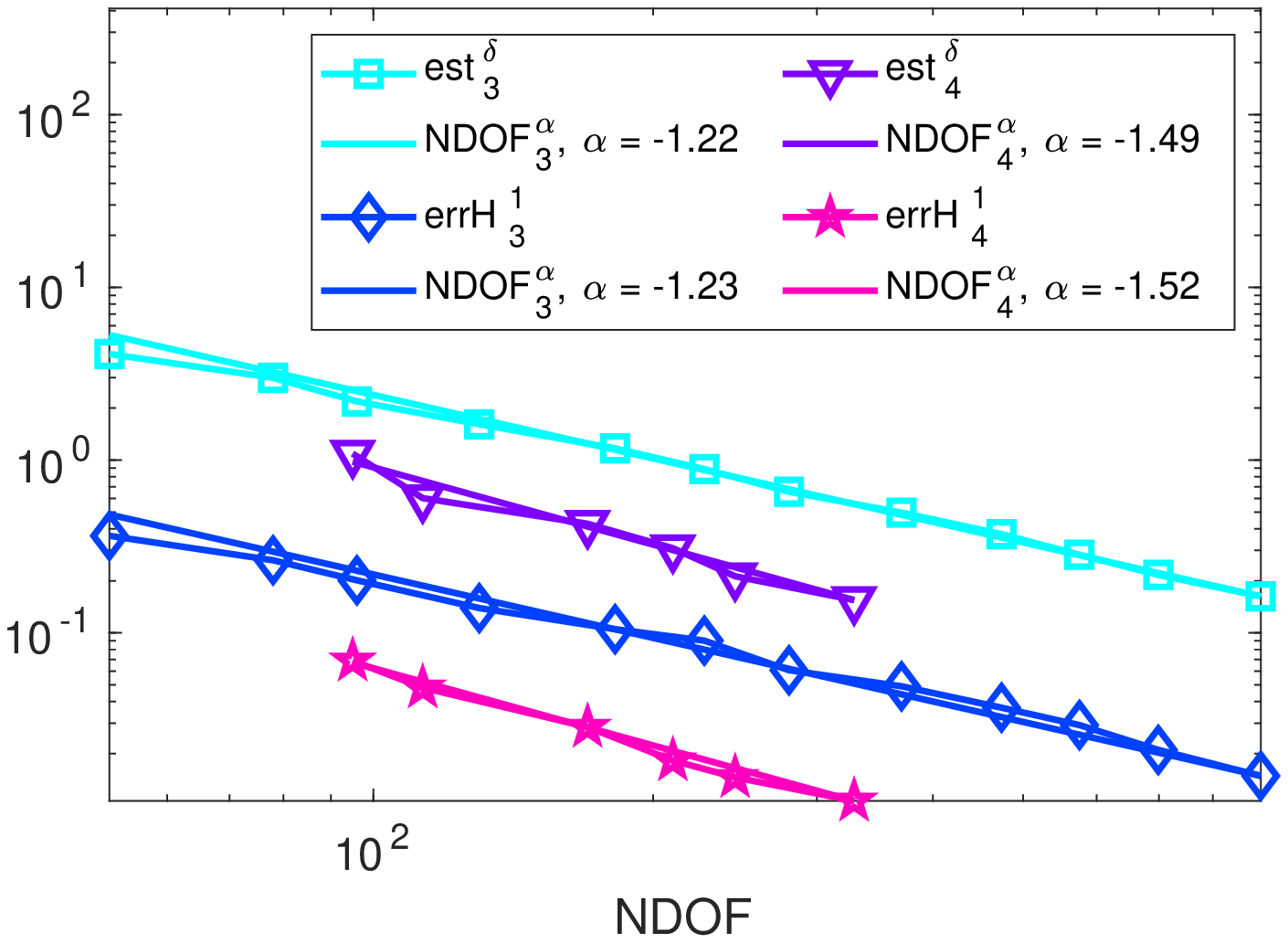}
    \label{fig:p2_TD}
  \end{subfigure}
  \caption{Problem 2: Rates of convergence ($\alpha$) for error and
    estimator with MaxMom refinement.}
  \label{fig:TipETut-Rates}
\end{figure}

\begin{table}[h]
  \centering \begin{tabular}{c|c|c|c|c}
TypeRef & Order 1 & Order 2 & Order 3 & Order 4 \\ 
\hline 
MaxMom & $\alpha = -0.48$ & $\alpha = -0.93$ & $\alpha = -1.22$ & $\alpha = -1.49$ \\ 
   & $\alpha_{err} = -0.49$ & $\alpha_{err} = -0.95$ & $\alpha_{err} = -1.23$ & $\alpha_{err} = -1.52$ \\ 
TrDir & $\alpha = -0.49$ & $\alpha = -0.95$ & $\alpha = -1.14$ & $\alpha = -1.49$ \\ 
   & $\alpha_{err} = -0.49$ & $\alpha_{err} = -0.92$ & $\alpha_{err} = -1.20$ & $\alpha_{err} = -1.52$ \\ 
MaxEdg & $\alpha = -0.49$ & $\alpha = -0.92$ & $\alpha = -1.32$ & $\alpha = -1.47$ \\ 
   & $\alpha_{err} = -0.47$ & $\alpha_{err} = -0.94$ & $\alpha_{err} = -1.35$ & $\alpha_{err} = -1.68$ \\ 
\end{tabular}

  \caption{Problem 2: Rates of convergence for the estimator ($\alpha$) and the error ($\alpha_{err}$) with refinement criteria MaxMom, TrDir and MaxEdg.}
  \label{tab:TipETut-Rates2}
\end{table}
\begin{table}
  \centering
  \begin{tabular}{c|c|c|c|c}
Step & $\#$ Cells & NDOF & $\varepsilon$ & PCG-It \\ 
\hline 
1 & 12 & 2 & 0.0296 & 1 \\ 
5 & 30 & 16 & 0.0472 & 6 \\ 
10 & 114 & 87 & 0.0541 & 13 \\ 
15 & 455 & 432 & 0.0552 & 25 \\ 
20 & 1875 & 1979 & 0.0557 & 55 \\ 
25 & 7752 & 8622 & 0.0552 & 114 \\ 
30 & 32221 & 36261 & 0.0531 & 234 \\ 
32 & 57076 & 64484 & 0.0520 & 309 \\ 
\end{tabular}
  \caption{Problem 2: MaxMom refinement Order 1.}
  \label{tab:TipETut-MM}
\end{table}
\begin{table}
  \centering
  \begin{tabular}{c|c|c|c|c}
Step & $\#$ Cells & NDOF & $\varepsilon$ & PCG-It \\ 
\hline 
1 & 12 & 2 & 0.0296 & 1 \\ 
5 & 30 & 16 & 0.0441 & 5 \\ 
10 & 113 & 92 & 0.0500 & 13 \\ 
15 & 433 & 428 & 0.0531 & 32 \\ 
20 & 1741 & 1933 & 0.0596 & 77 \\ 
25 & 7047 & 8221 & 0.0581 & 152 \\ 
30 & 28933 & 34475 & 0.0574 & 411 \\ 
33 & 67557 & 80552 & 0.0562 & 542 \\ 
\end{tabular}
  \caption{Problem 2: TrDir refinement Order 1.}
  \label{tab:TipETut-TD}
\end{table}
\begin{table}
  \centering
  \begin{tabular}{c|c|c|c|c}
Step & $\#$ Cells & NDOF & $\varepsilon$ & PCG-It \\ 
\hline 
1 & 12 & 2 & 0.0296 & 1 \\ 
5 & 30 & 16 & 0.0476 & 6 \\ 
10 & 113 & 89 & 0.0555 & 12 \\ 
15 & 450 & 418 & 0.0578 & 24 \\ 
20 & 1847 & 1885 & 0.0614 & 55 \\ 
25 & 7737 & 8203 & 0.0636 & 114 \\ 
30 & 32385 & 34978 & 0.0679 & 243 \\ 
31 & 43231 & 47014 & 0.0681 & 286 \\ 
\end{tabular}
  \caption{Problem 2: MaxEdg refinement Order 1.}
  \label{tab:TipETut-ME}
\end{table}


The geometry and the parameters of this test problem are described in
detail in \cite[Section 6.1]{BBBPS}.  The DFN is composed by three
fractures and three traces, defined as

\begin{center}
\begin{tabular}{p{.38\linewidth}|p{.49\linewidth}}
{\begin{align*}%
     F_1 &= \left(-1,\frac12\right) \!\times\! \left(-1,1\right) \!\times\! \{0\} \,,%
     \\%
     F_2 &= (-1,0) \!\times\! \{0\} \!\times\! (-1,1) \,,%
     \\%
     F_3 &= \left\{\frac12\right\} \!\times\! (-1,1) \!\times\! (-1,1) \,,%
 \end{align*}}
         &
           {\begin{align*} \Gamma_1 &= F_1 \cap F_2 = \left(-1,\frac12\right) \!\times\! \{0\} \!\times\!
                                      \{0\} \,,%
              \\%
              \Gamma_2 &= F_1 \cap F_3 = \left(-1,0\right) \!\times\! \{0\} \!\times\!
                         \left(-1,1\right) \,,%
              \\%
              \Gamma_3 &= F_2 \cap F_3 = \left\{-\frac12\right\} \!\times\! \{0\} \!\times\!
                         (-1,1)\,.\end{align*}}%
\end{tabular}
\end{center}
The chosen exact solution is defined on each fracture as follows:
\begin{align*}
  h_1(x,y,z) &= -\frac{1}{10} \left( x+\frac12 \right)
               \left( 8 x y \left( x^2+y^2 \right)
               \mathrm{arctan2}(y,x) + x^3 \right) \,,
  \\
  h_2(x,y,z) &= -\frac{1}{10} \left( x + \frac12 \right) x^3
               \left( 1 - 8 \pi \abs{z} \right) \,,
  \\
  h_3(x,y,z) &= y(y-1)(y+1)(z-1)z \,,
\end{align*}
and all fractures have transmissivities set to $1$.

In Figure \ref{fig:TipETut_DFN-mesh}, we report some of the meshes
generated by the refinement process. Similarly to the previous test, a
refinement is induced also far from traces due to the presence of a
non-null forcing term.

In Figure \ref{fig:TipETut-Rates} we report the rates of convergence ($\alpha$)
of the error and of the error estimator considering the
strategy MaxMom. The rates of convergence for the strategies MaxMom, TrDir and MaxEdg are shown
in the Table~\ref{tab:TipETut-Rates2}. The convergence rates are computed on the
basis of the last five refinement iterations: again we can remark a
very good agreement between the error and the estimator. The sub-optimal rates of convergence with higher VEM orders is still due to the bounded regularity of the solution.

In Tables \ref{tab:TipETut-MM}-\ref{tab:TipETut-ME} we report the same quantities reported in Tables
\ref{tab:TipEnrich-MM}-\ref{tab:TipEnrich-ME} in order to describe
the refinement process. We observe again that the effectivity index is
almost independent of the meshsize since, as the mesh starts to have a
sufficient number of DOFs, it displays small variations.

 
\section{Numerical Results on a realistic DFN}
\label{sec:num_res_DFN}
\begin{figure}[h!]
  \centering \includegraphics[width=\linewidth]{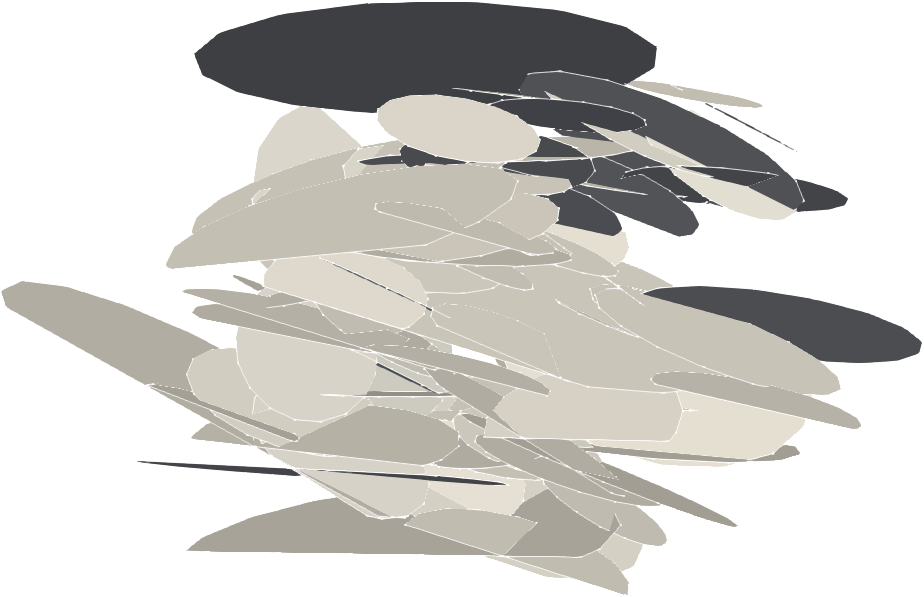}
  \caption{DFN with 86 fractures.}
  \label{fig:dfn86}
\end{figure}

In this section we discuss the results obtained by the four presented
refinement strategies when applied to more realistic DFNs.  The
geometry of the considered DFNs is fixed and is composed by 86
fractures and 159 traces, with a maximum number of traces per fracture
equal to 11 and a mean value of traces per fracture equal to 1.85, see
Figure \ref{fig:dfn86}. We consider two test cases where the
transmissivities of the fractures are sampled from two log-normal
distributions having standard deviations equal to $10$ and $10^4$,
respectively. These two problems are tagged with the labels DFN86E01
and DFN86E04.

The problems considered have no forcing terms, and the flux is driven
by the presence of two Dirichlet boundary conditions (10 on the boundaries at $x=0$ and 0 at the boundaries at $x=1000$), whereas homogeneous Neumann boundary conditions are imposed on all the other boundaries.

In the following analysis we display the convergence rates of the
error estimate with respect to the number of degrees of freedom and we
assess how the refinement strategies impact on the aspect ratio of the
cells of a selected fracture (Fracture $72$), and on the iterations of the
preconditioned conjugate gradient method used to solve the linear
system.

\begin{figure}[h!]
  \centering
  \begin{subfigure}[b]{0.49\linewidth}
    \includegraphics[width=\linewidth]{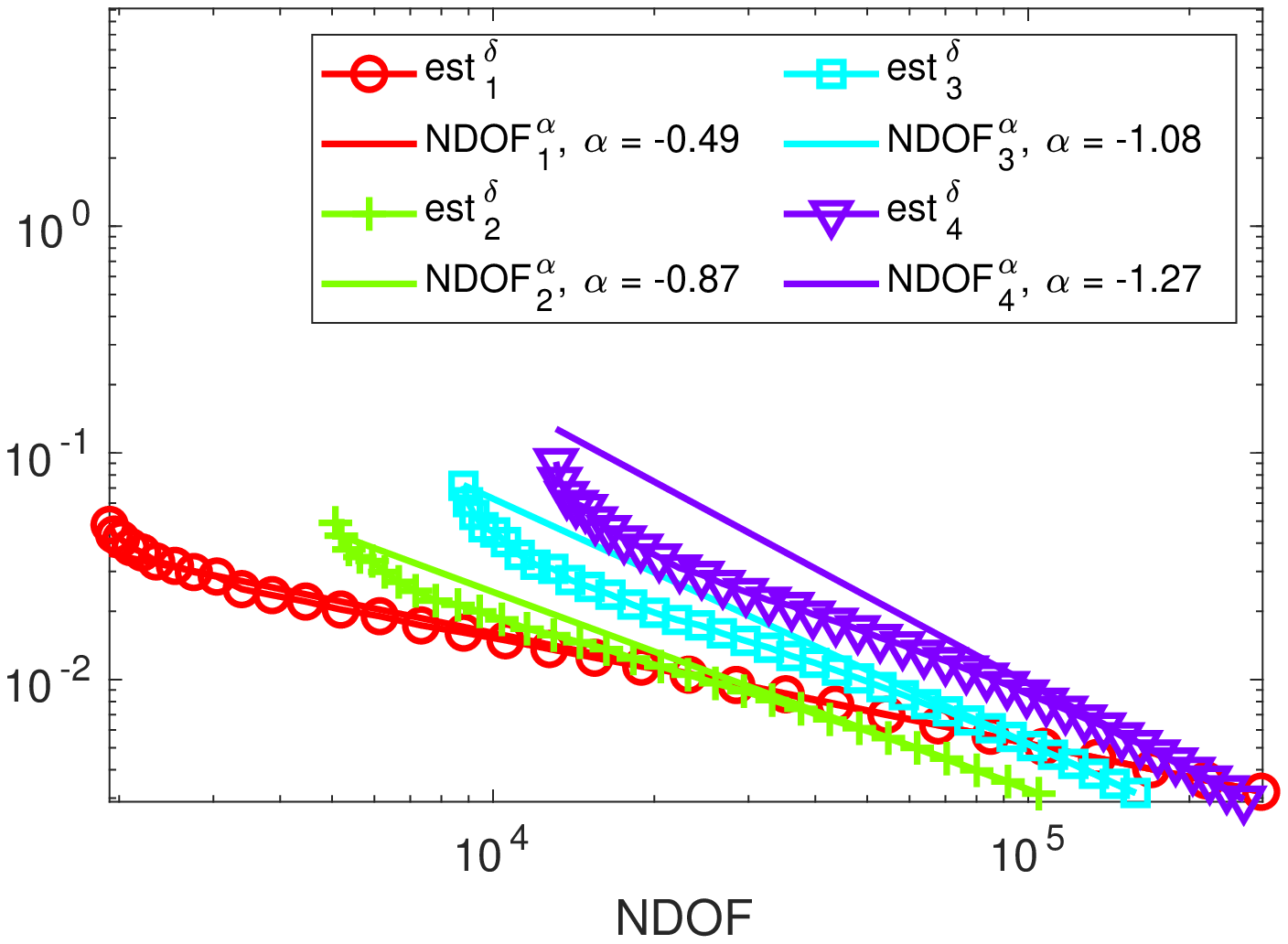}
    \caption{MaxMom refinement.}
    \label{fig:ten1_estDof0}
  \end{subfigure}
  \begin{subfigure}[b]{0.49\linewidth}
    \includegraphics[width=\linewidth]{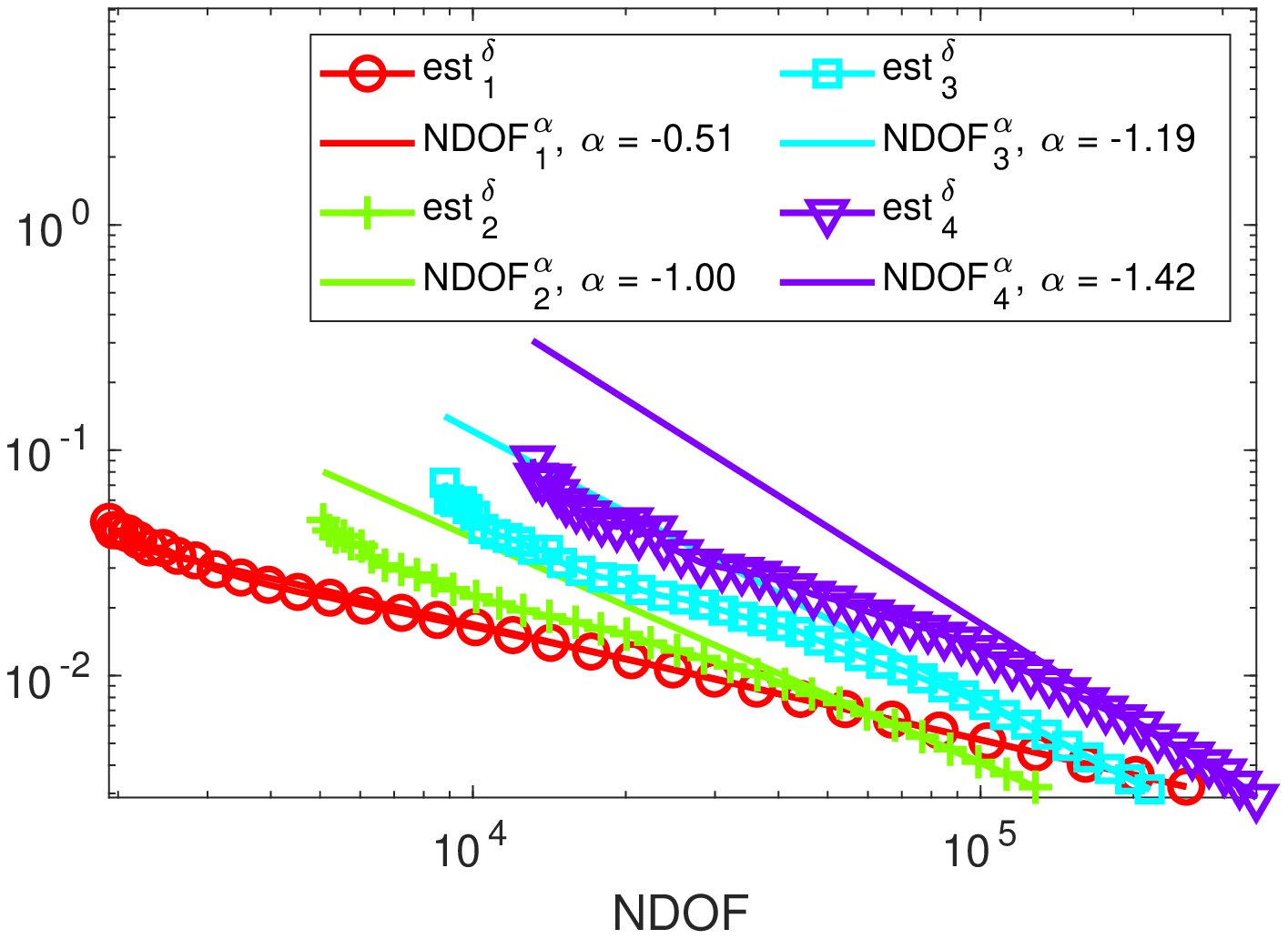}
    \caption{TrDir refinement.}
    \label{fig:ten1_estDof1}
  \end{subfigure}
  \begin{subfigure}[b]{0.49\linewidth}
    \includegraphics[width=\linewidth]{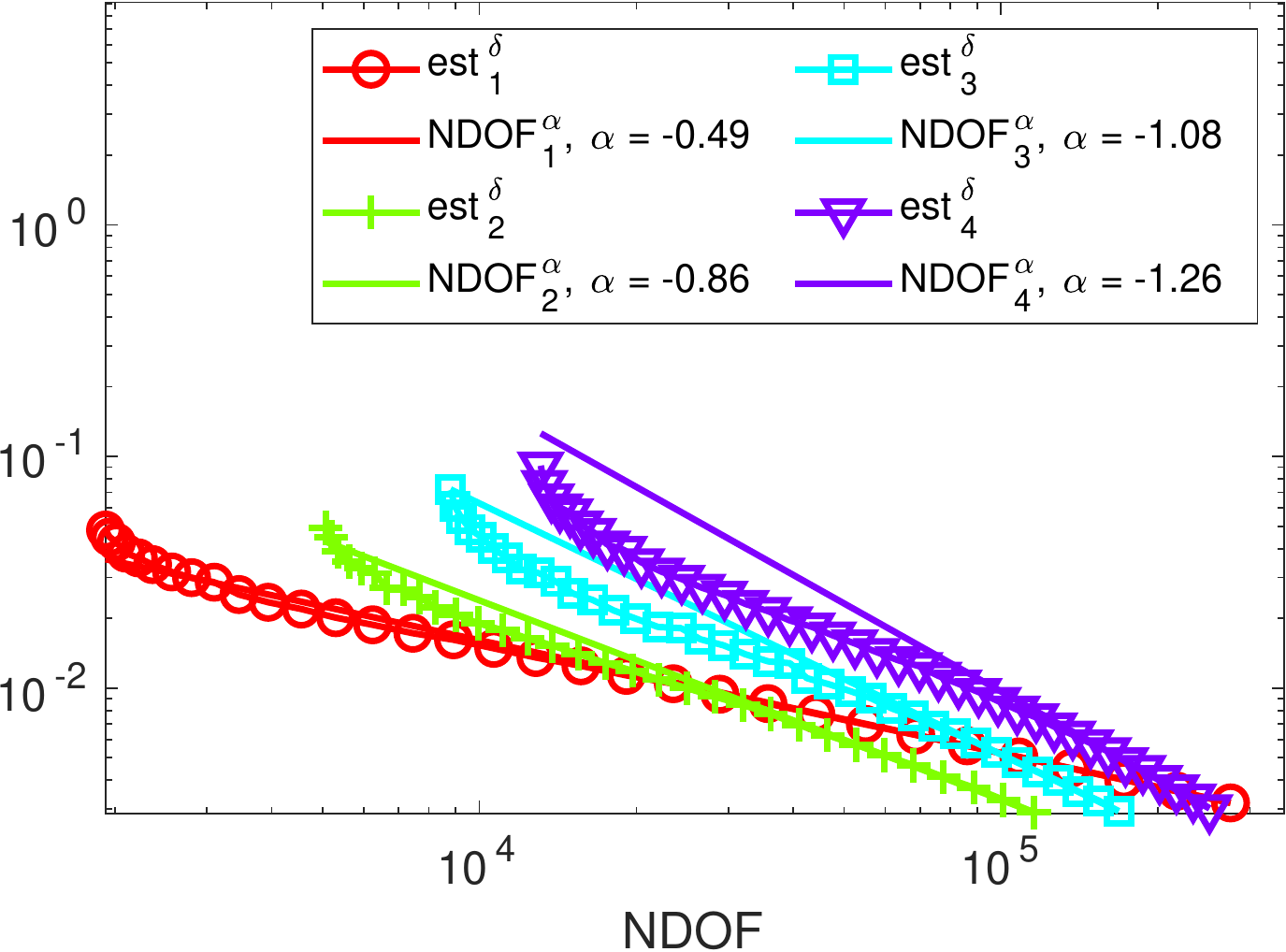}
    \caption{MaxPnt refinement.}
    \label{fig:ten1_estDof2}
  \end{subfigure}
  \begin{subfigure}[b]{0.49\linewidth}
    \includegraphics[width=\linewidth]{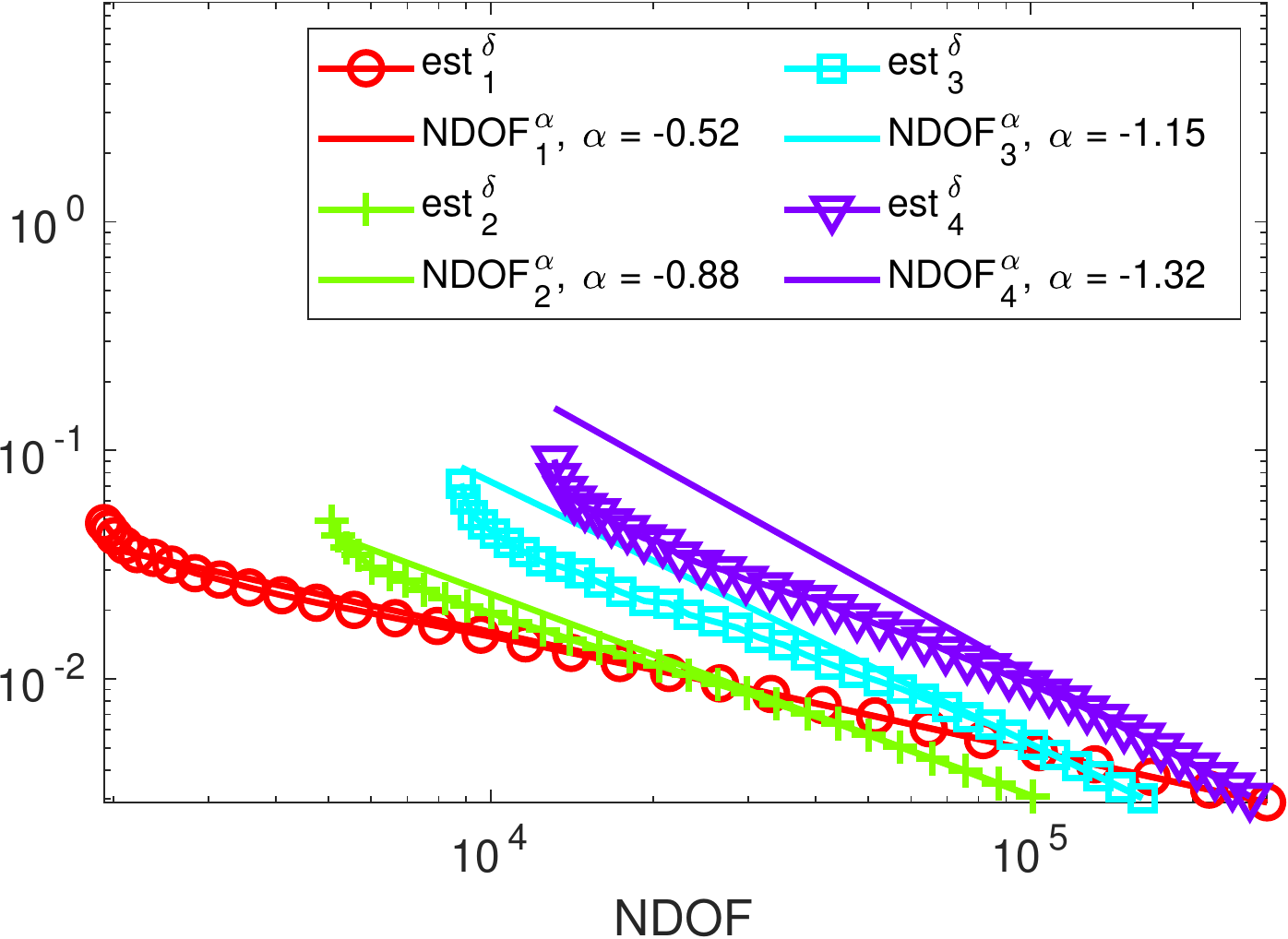}
    \caption{MaxEdg refinement.}
    \label{fig:ten1_estdDof3}
  \end{subfigure}
  \caption{DFN86E01: Estimator vs. NDOF}
  \label{fig:ten1_estdDof}
\end{figure}

\begin{figure}[h!]
  \centering
  \begin{subfigure}[b]{0.49\linewidth}
    \includegraphics[width=\linewidth]{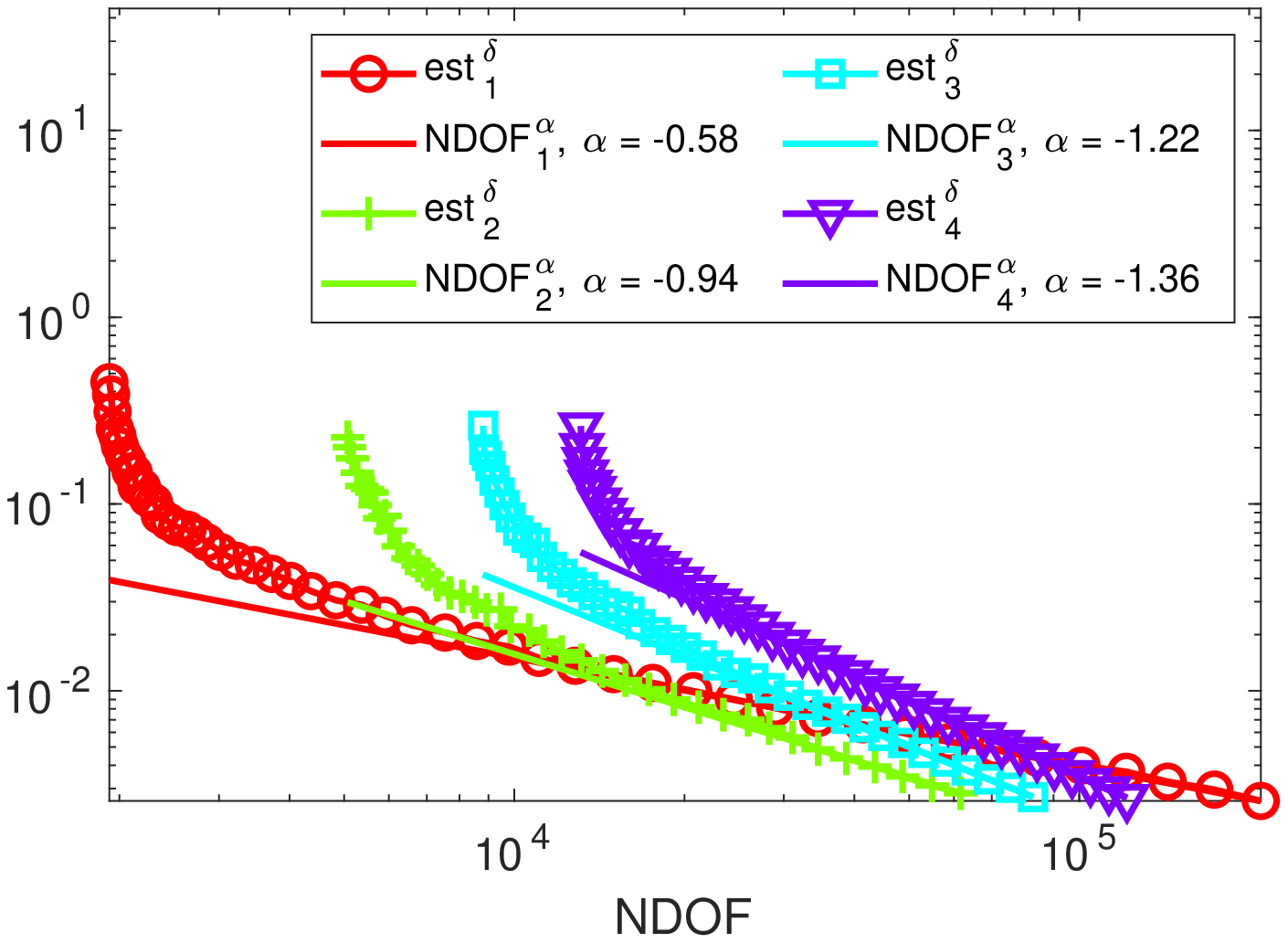}
    \caption{MaxMom refinement.}
    \label{fig:ten4_estDof0}
  \end{subfigure}
  \begin{subfigure}[b]{0.49\linewidth}
    \includegraphics[width=\linewidth]{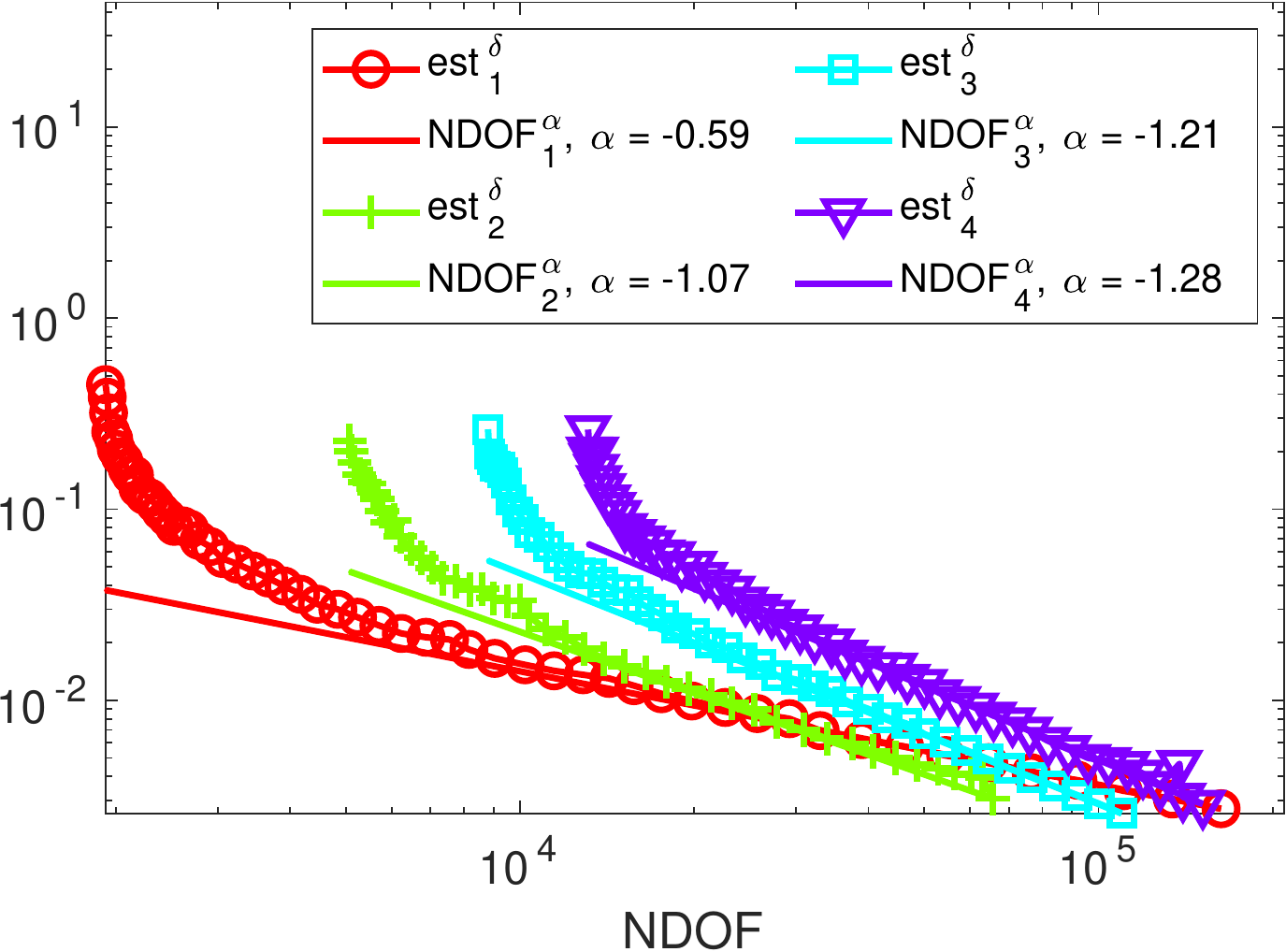}
    \caption{TrDir refinement.}
    \label{fig:ten4_estDof1}
  \end{subfigure}
  \begin{subfigure}[b]{0.49\linewidth}
    \includegraphics[width=\linewidth]{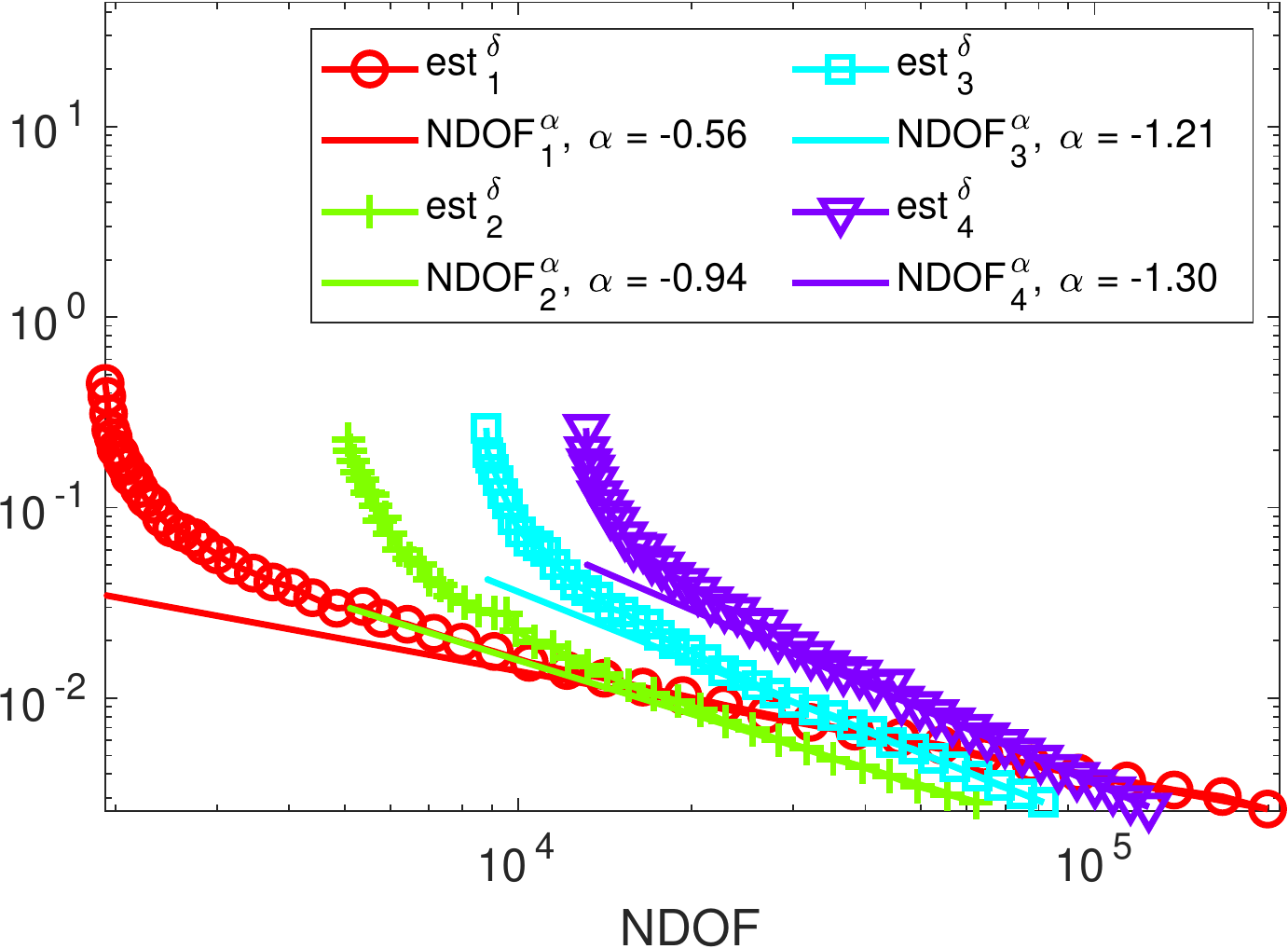}
    \caption{MaxPnt refinement.}
    \label{fig:ten4_estDof2}
  \end{subfigure}
  \begin{subfigure}[b]{0.49\linewidth}
    \includegraphics[width=\linewidth]{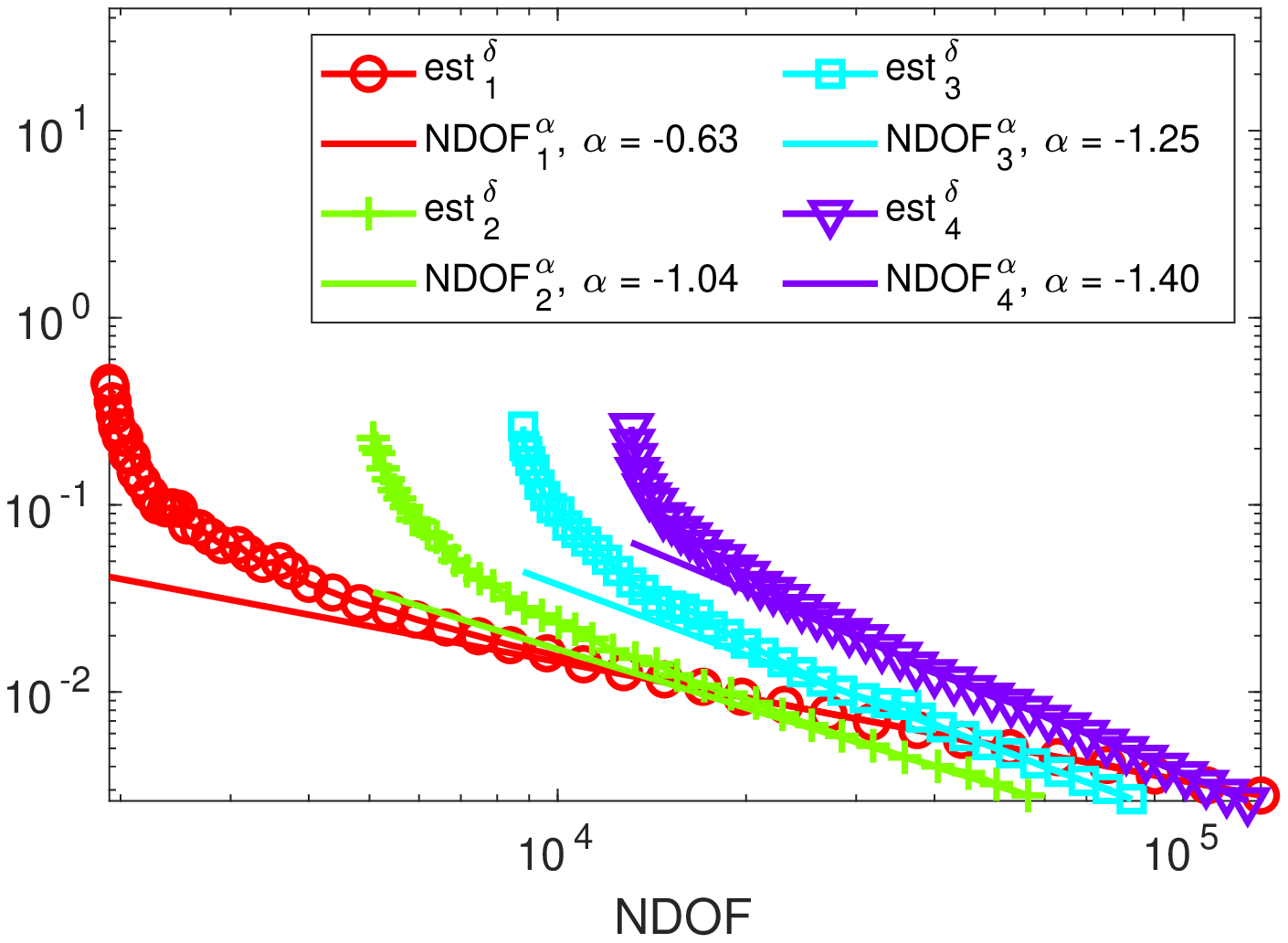}
    \caption{MaxEdg refinement.}
    \label{fig:ten4_estDof3}
  \end{subfigure}
  \caption{DFN86E04: Estimator vs. NDOF}
  \label{fig:ten4_estdDof}
\end{figure}

In Figures \ref{fig:ten1_estdDof} and \ref{fig:ten4_estdDof} we
display the behaviour of the estimators with respect to the number of
DOFs for the four refinement strategies considered and we report the
slope of the estimator for each VEM order, for the two problems. We
can observe that all the refinement strategies display an optimal
asymptotic rate of convergence up to VEM order 2 ($-0.5$ for the VEM
Order 1 and $-1.0$ for the VEM Order 2). For higher VEM orders the
bounded regularity constraints the rates of convergence. The trend of
the rate of convergence results similar even if the fluxes in the
these two DFNs are completely different.

\begin{figure}
  \begin{subfigure}[b]{0.49\linewidth}
  \centering
  \includegraphics[width=\linewidth]{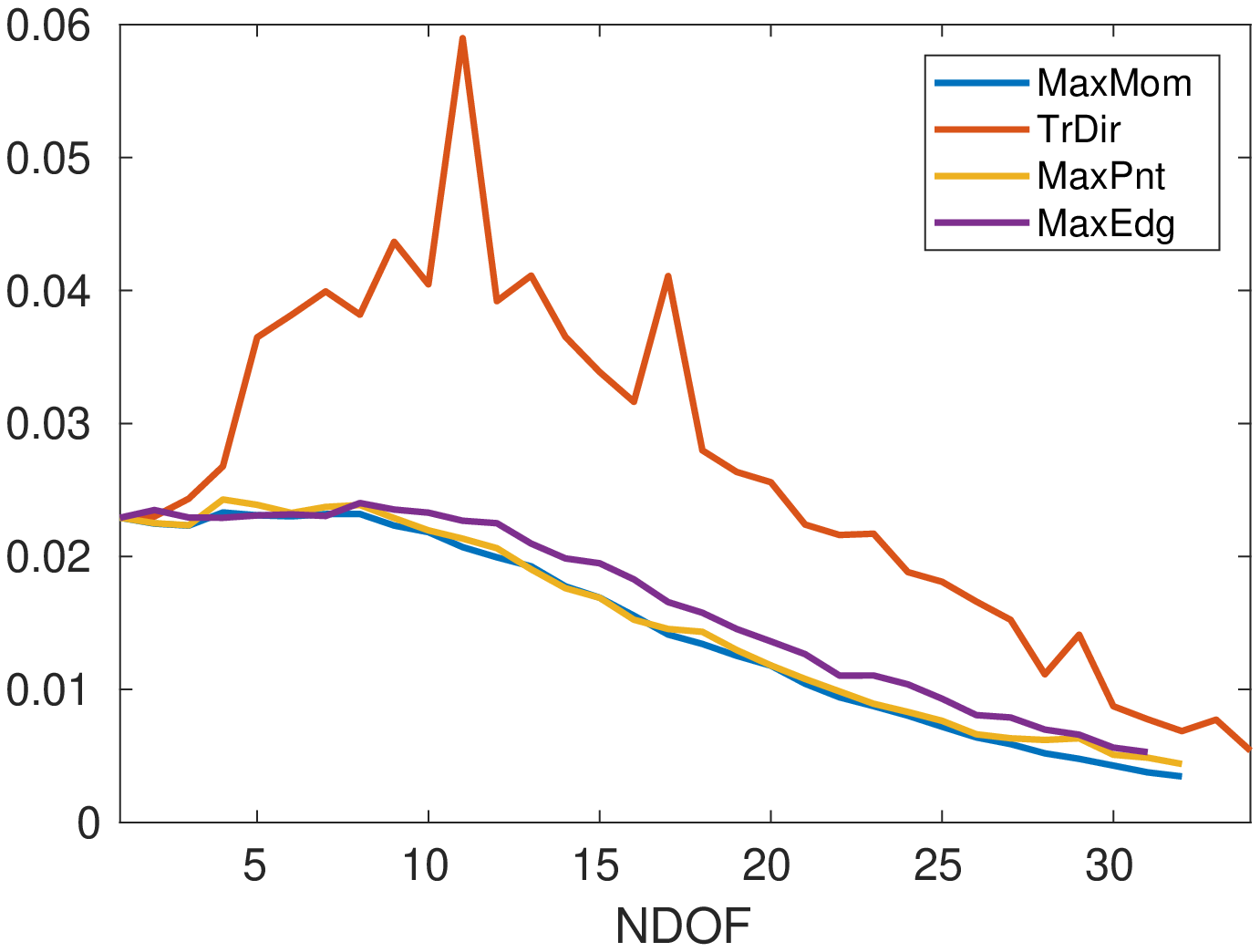}
  \caption{DFN86E01}
  \label{fig:ten1_it}
  \end{subfigure}
  \begin{subfigure}[b]{0.49\linewidth}
  \centering
  \includegraphics[width=\linewidth]{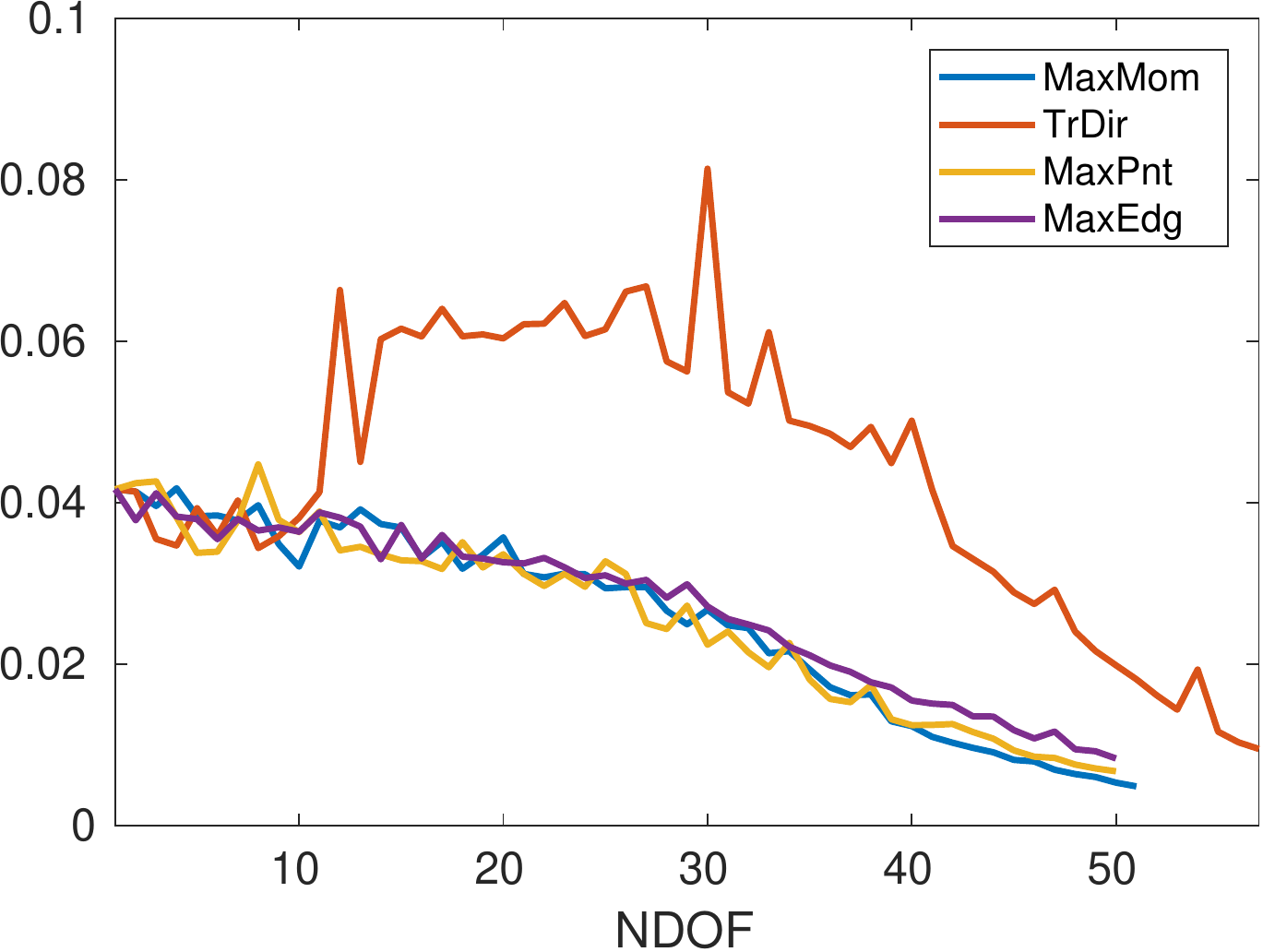}
  \caption{DFN86E04}
  \label{fig:ten4_it}
\end{subfigure}
\caption{Ratio between PCG-It and NDOF}
\label{fig:it}
\end{figure}

In Figure \ref{fig:it} we plot the
behaviour of the ratio PCG-It/NDOF. After the initial noisy behaviour
we can observe a decreasing asymptotical trend for all the considered
refinement strategies. These plots highlight the advantages of a
suitably refined mesh also on the performances of the linear solver.

\begin{figure}[h!]
  \centering \includegraphics[width=0.5\linewidth]{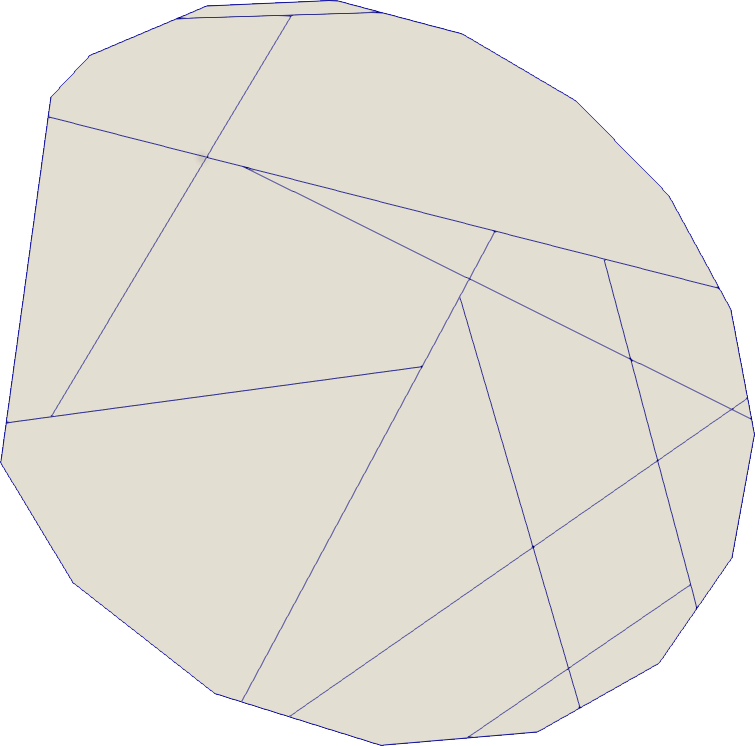}
  \caption{DFN86E01 - DFN86E04: initial mesh on Fracture 72.}
  \label{fig:ten1_f72s0}
\end{figure}

\begin{figure}[h!]
  \centering
  \begin{subfigure}[b]{0.4\linewidth}
    \includegraphics[width=\linewidth]{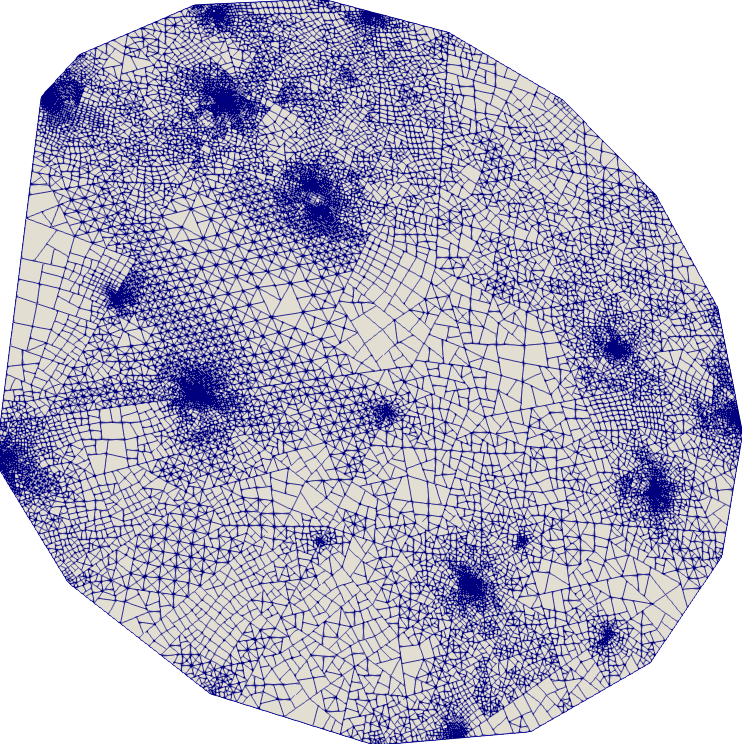}
    \caption{MaxMom refinement.}
  \end{subfigure}
  \begin{subfigure}[b]{0.4\linewidth}
    \includegraphics[width=\linewidth]{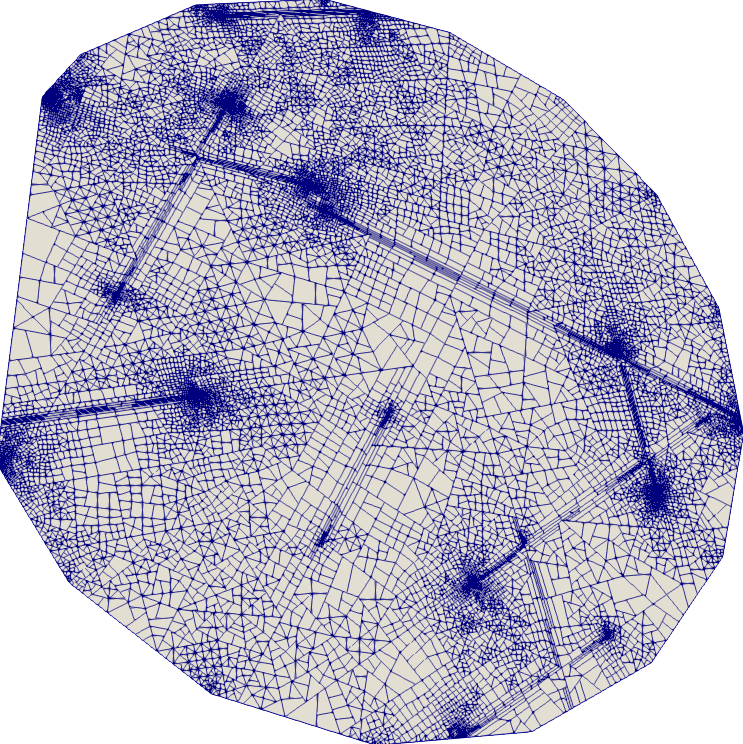}
    \caption{TrDir refinement.}
  \end{subfigure}
  \begin{subfigure}[b]{0.4\linewidth}
    \includegraphics[width=\linewidth]{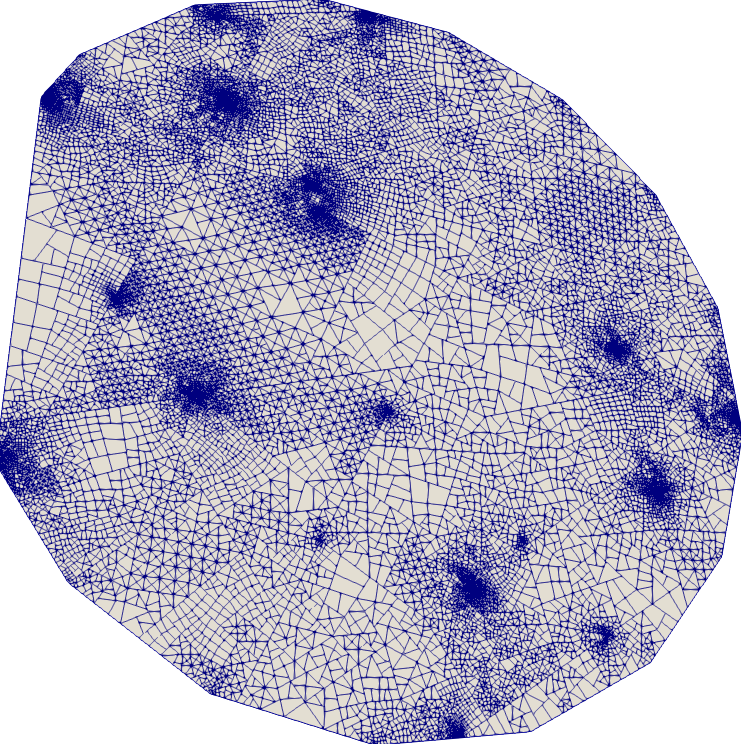}
    \caption{MaxPnt refinement.}
  \end{subfigure}
  \begin{subfigure}[b]{0.4\linewidth}
    \includegraphics[width=\linewidth]{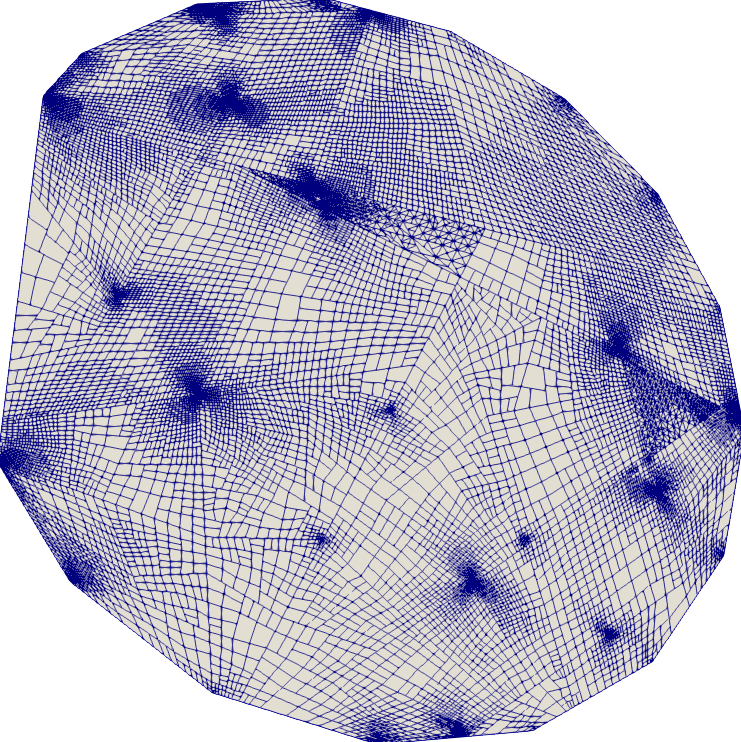}
    \caption{MaxEdg refinement.}
  \end{subfigure}
  \caption{DFN86E01: final mesh on Fracture 72.}
  \label{fig:ten1f72slast}
\end{figure}

\begin{figure}[h!]
  \centering
  \begin{subfigure}[b]{0.4\linewidth}
    \includegraphics[width=\linewidth]{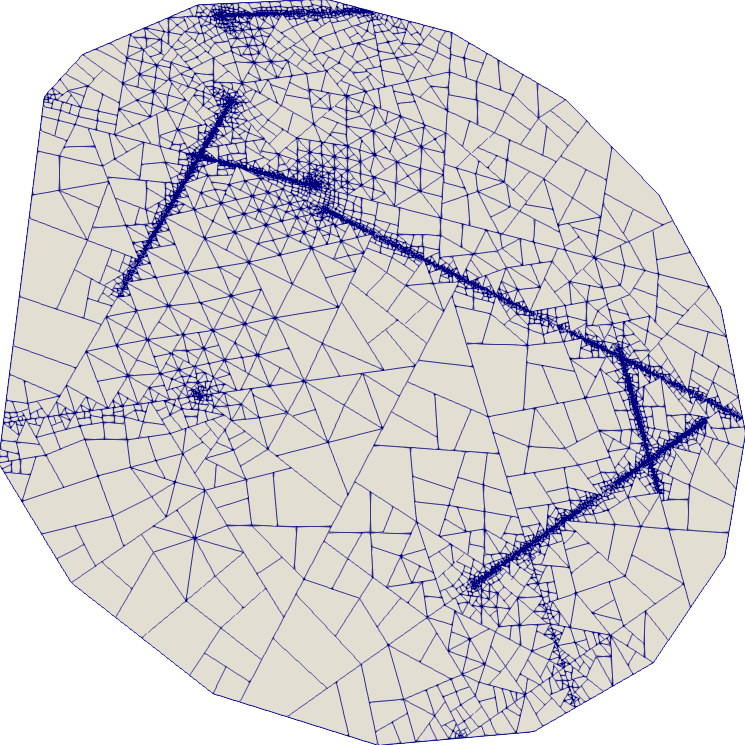}
    \caption{MaxMom refinement.}
  \end{subfigure}
  \begin{subfigure}[b]{0.4\linewidth}
    \includegraphics[width=\linewidth]{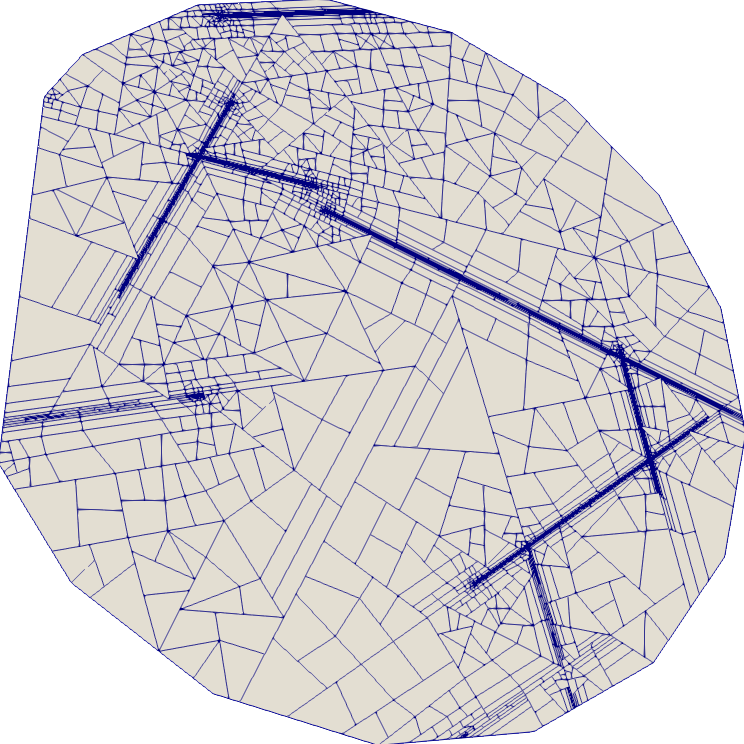}
    \caption{TrDir refinement.}
  \end{subfigure}
  \begin{subfigure}[b]{0.4\linewidth}
    \includegraphics[width=\linewidth]{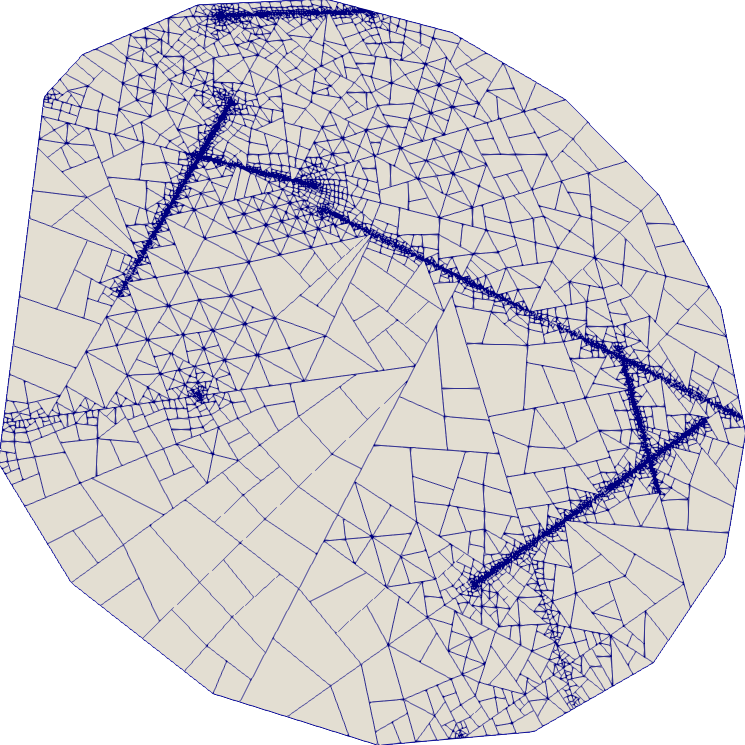}
    \caption{MaxPnt refinement.}
  \end{subfigure}
  \begin{subfigure}[b]{0.4\linewidth}
    \includegraphics[width=\linewidth]{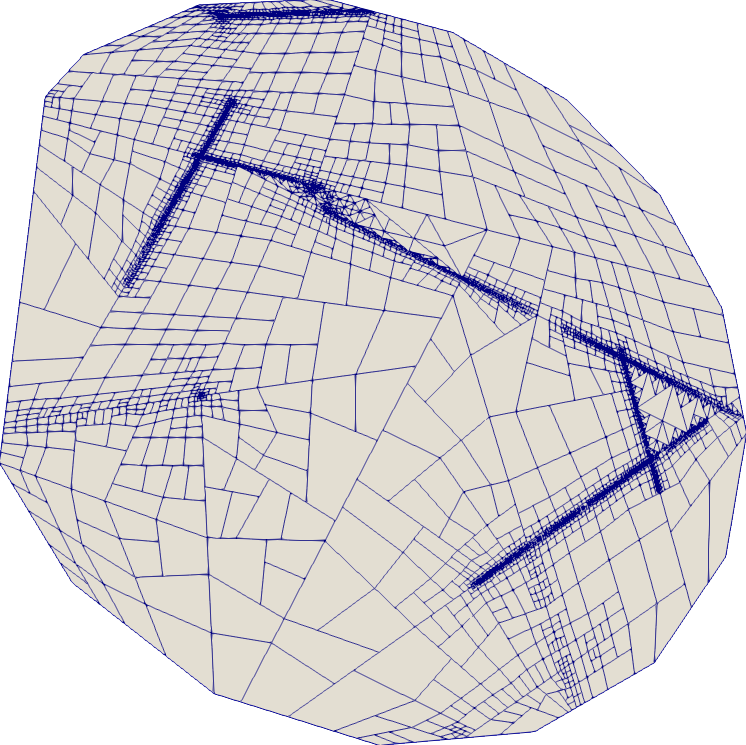}
    \caption{MaxEdg refinement.}
  \end{subfigure}
  \caption{DFN86E04: final mesh on Fracture 72.}
  \label{fig:ten4f72slast}
\end{figure}

\begin{figure}[h!]
  \centering
  \begin{subfigure}[b]{0.4\linewidth}
    \includegraphics[width=\linewidth]{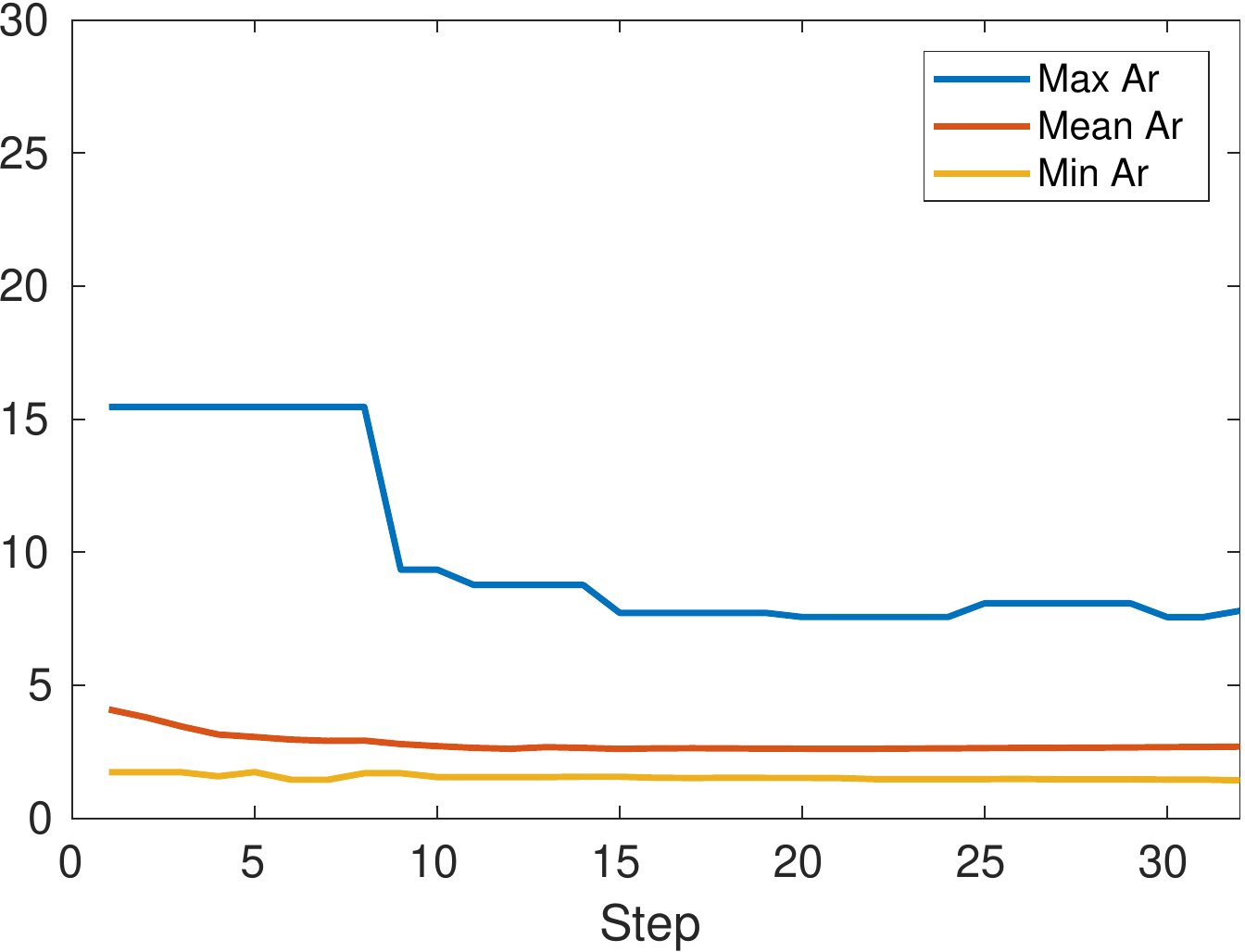}
    \caption{MaxMom refinement.}
  \end{subfigure}
  \begin{subfigure}[b]{0.4\linewidth}
    \includegraphics[width=\linewidth]{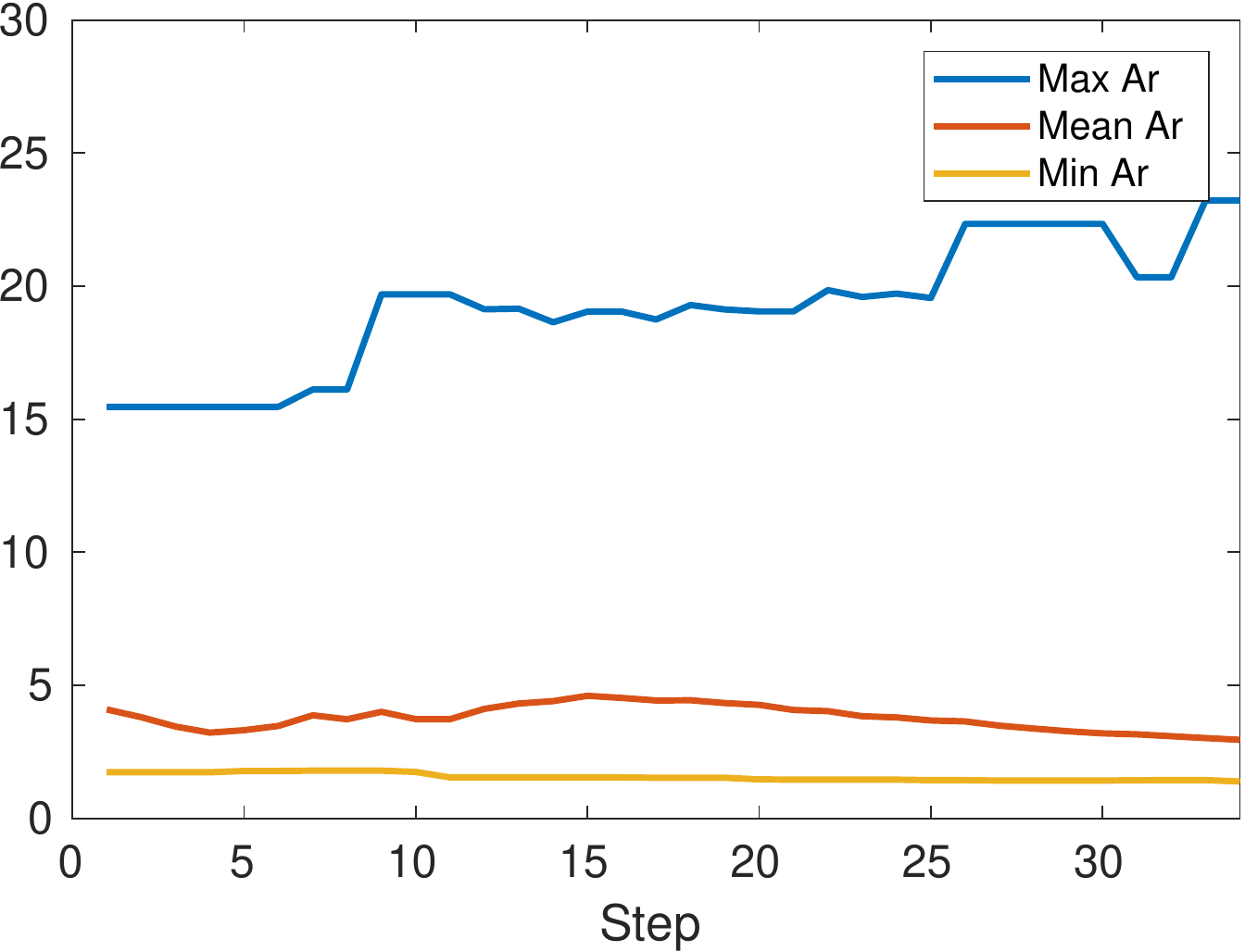}
    \caption{TrDir refinement.}
  \end{subfigure}
  \begin{subfigure}[b]{0.4\linewidth}
    \includegraphics[width=\linewidth]{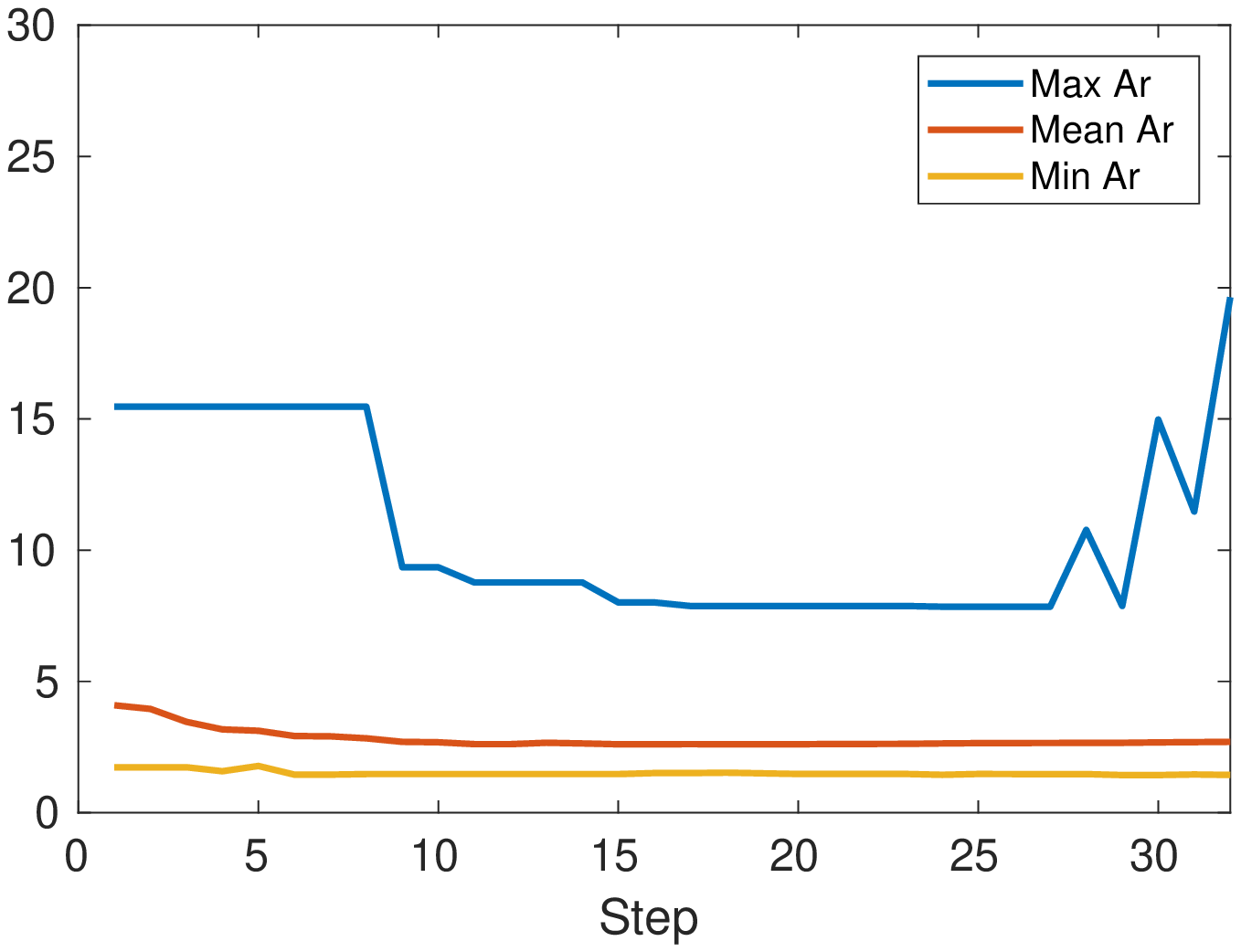}
    \caption{MaxPnt refinement.}
  \end{subfigure}
  \begin{subfigure}[b]{0.4\linewidth}
    \includegraphics[width=\linewidth]{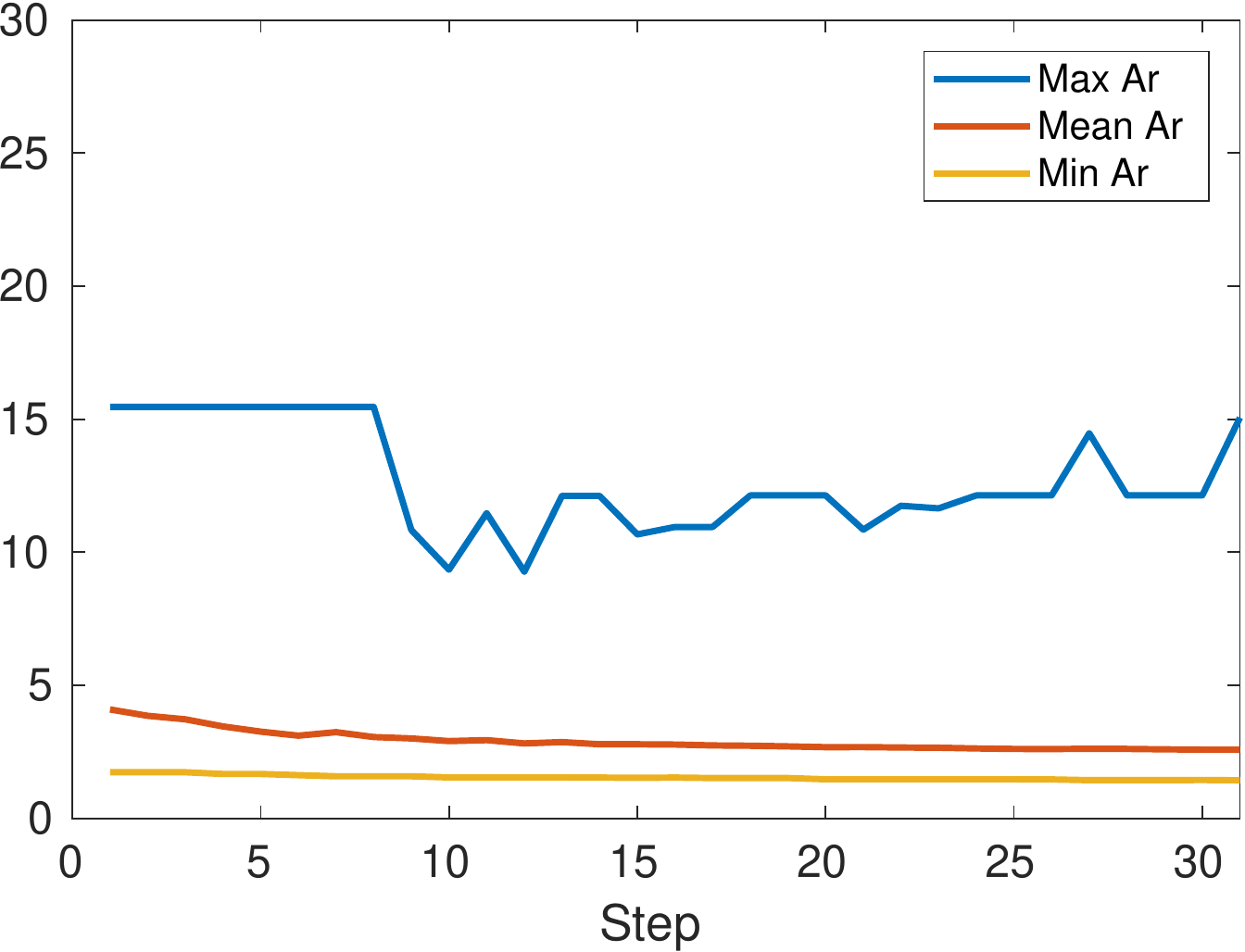}
    \caption{MaxEdg refinement.}
  \end{subfigure}
  \caption{DFN86E01: aspect ratio statistics of the mesh on Fracture
    72.}
  \label{fig:ten1f72sAr}
\end{figure}

\begin{figure}[h!]
  \centering
  \begin{subfigure}[b]{0.4\linewidth}
    \includegraphics[width=\linewidth]{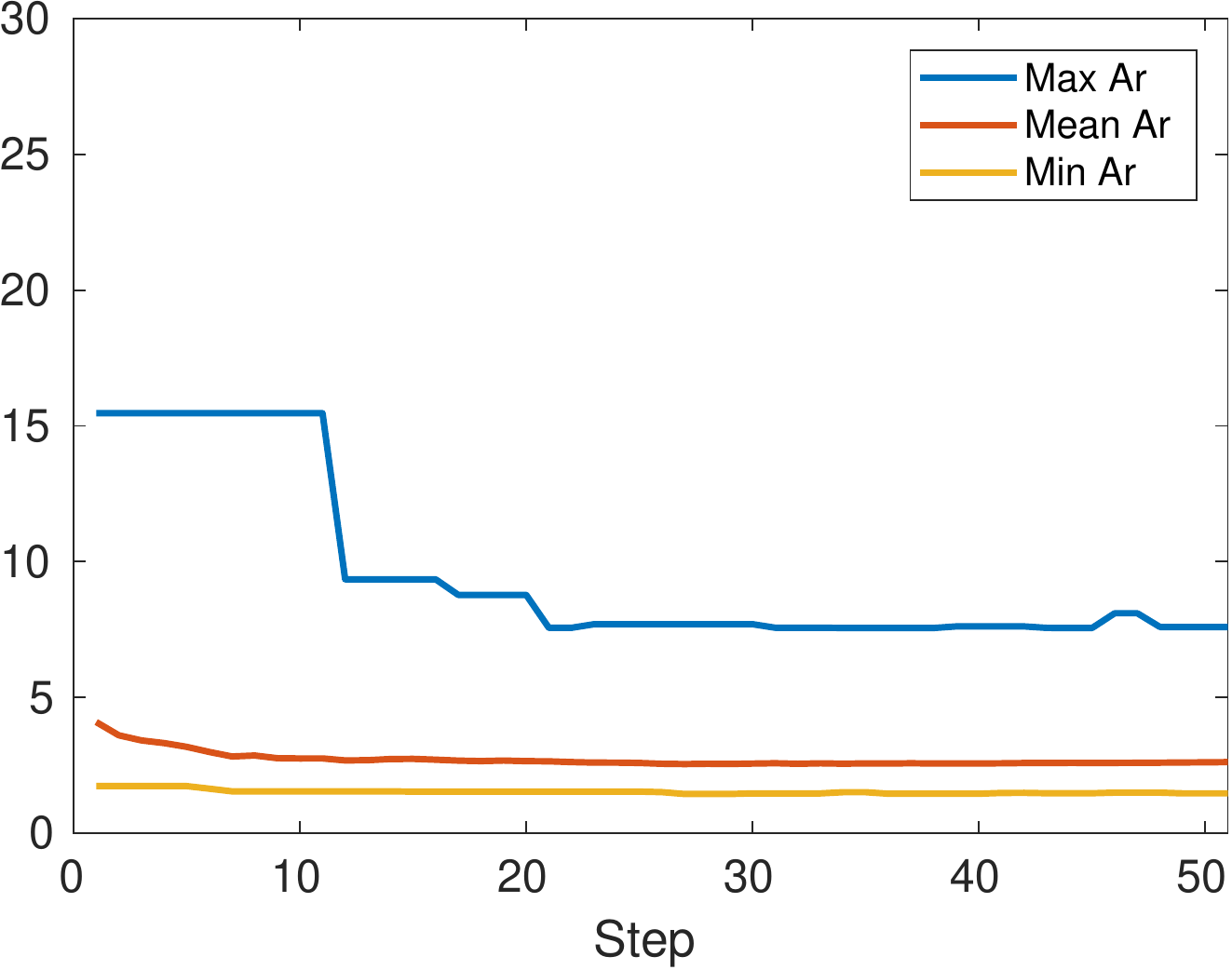}
    \caption{MaxMom refinement.}
  \end{subfigure}
  \begin{subfigure}[b]{0.4\linewidth}
    \includegraphics[width=\linewidth]{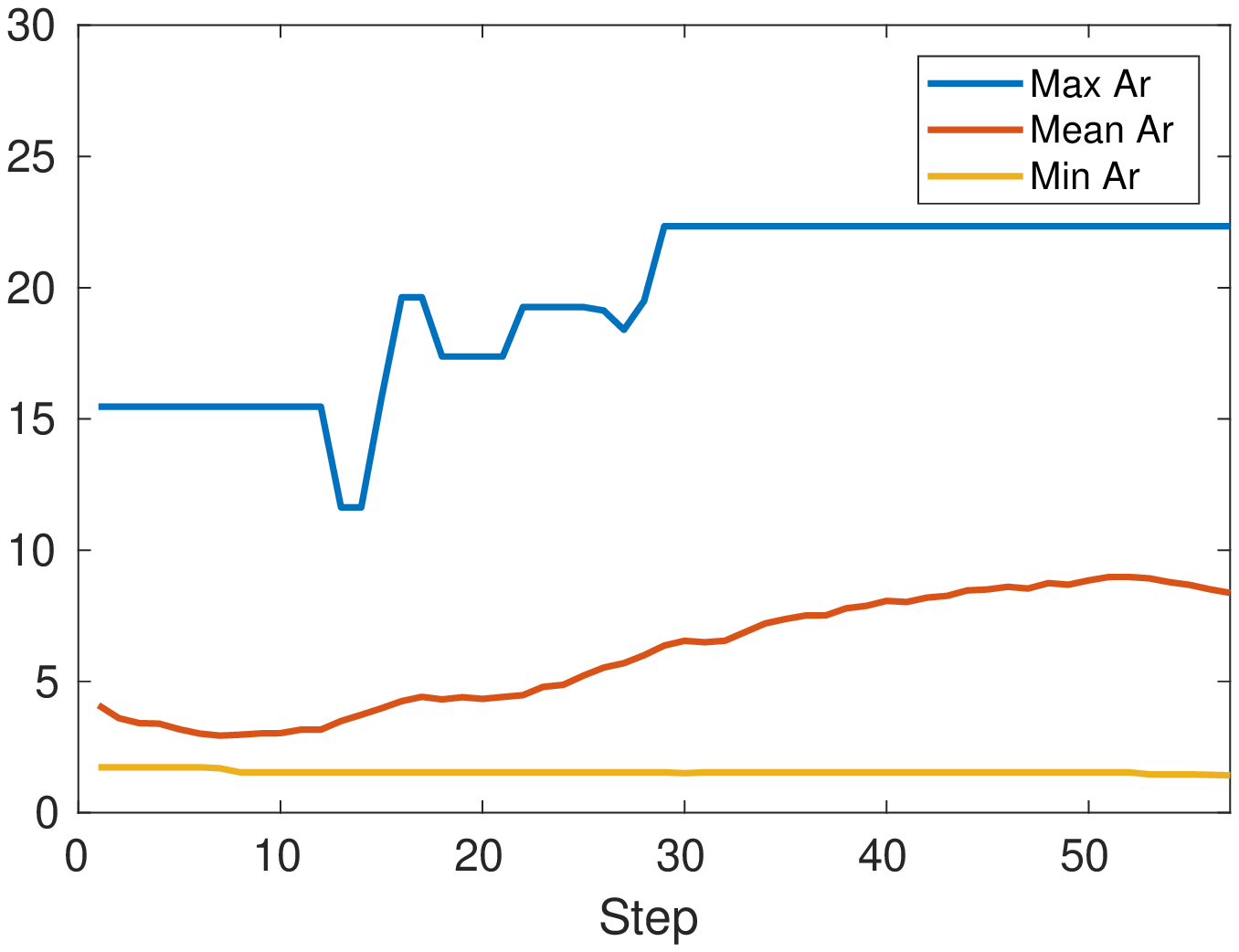}
    \caption{TrDir refinement.}
  \end{subfigure}
  \begin{subfigure}[b]{0.4\linewidth}
    \includegraphics[width=\linewidth]{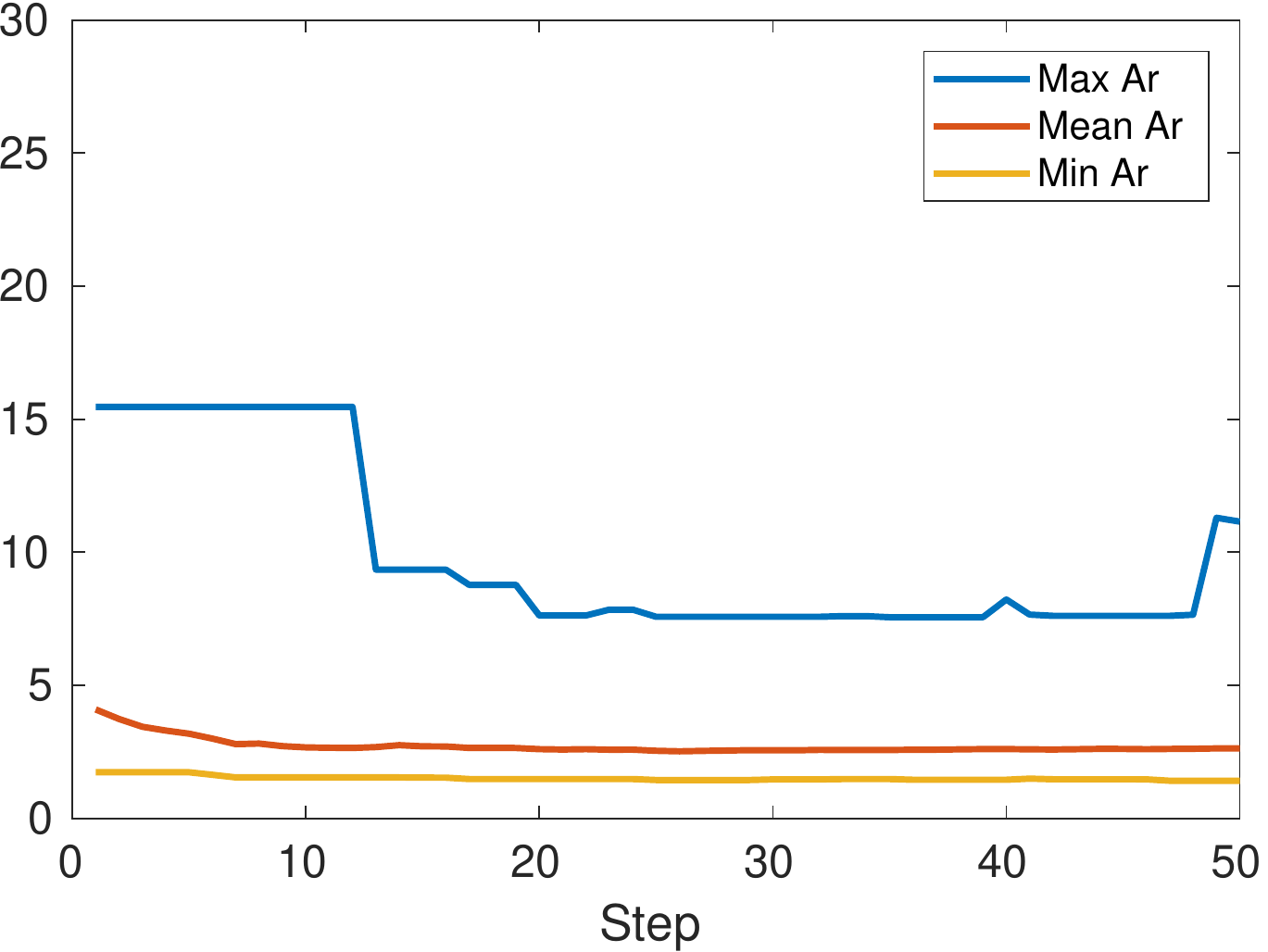}
    \caption{MaxPnt refinement.}
  \end{subfigure}
  \begin{subfigure}[b]{0.4\linewidth}
    \includegraphics[width=\linewidth]{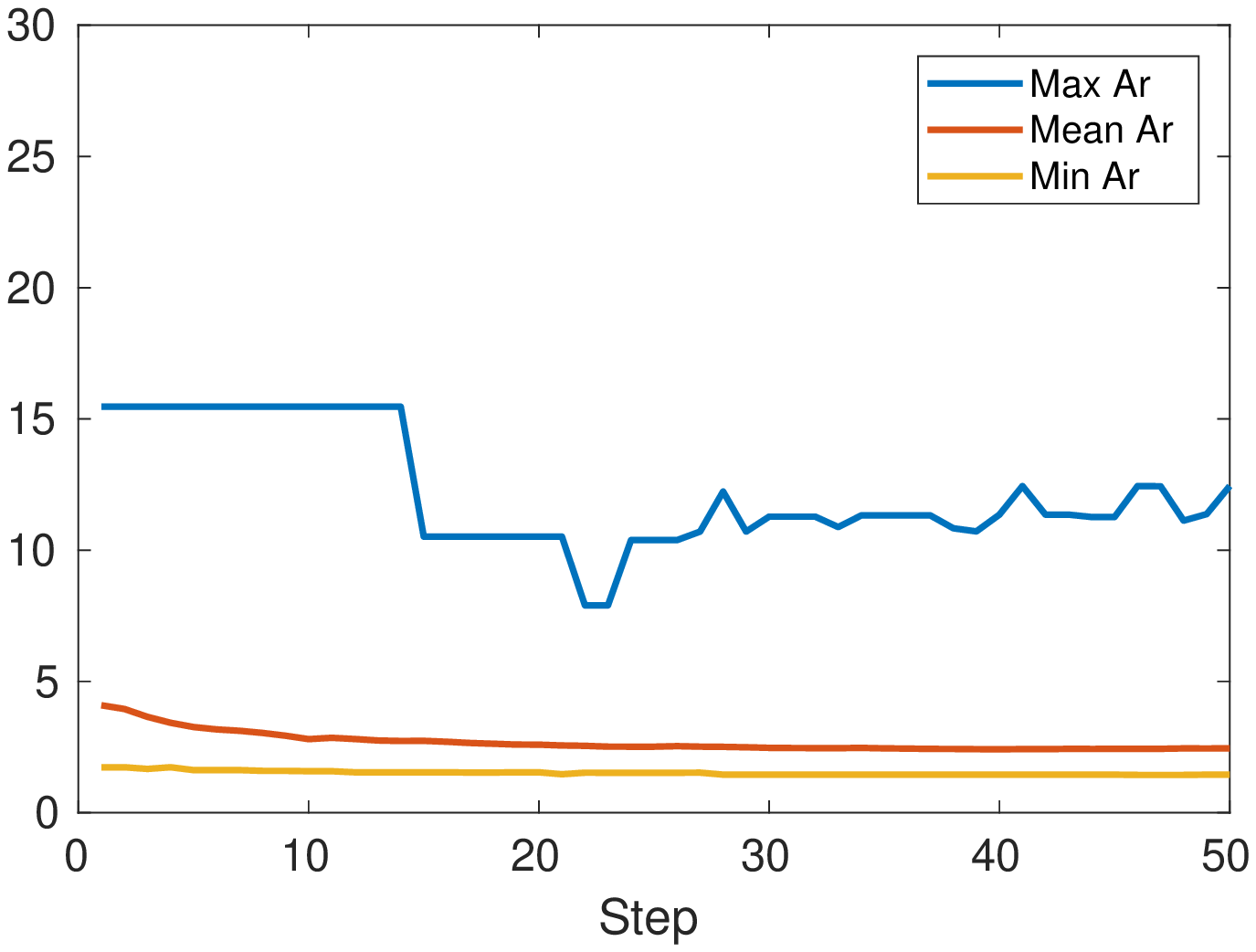}
    \caption{MaxEdg refinement.}
  \end{subfigure}
  \caption{DFN86E04: aspect ratio statistics of the mesh on Fracture
    72.}
  \label{fig:ten4f72Ar}
\end{figure}

In Figure~\ref{fig:ten1_f72s0} we report the minimal mesh on Fracture
72: this is the common initial mesh for all the
refinement strategies and for both DFN86E01 and DFN86E04. In Figures
\ref{fig:ten1f72slast} and \ref{fig:ten4f72slast} we display the
meshes produced by all the considered refinement strategies. The TrDir
strategy produces a stronger and sharper refinement along the
traces. The switch from TrDir to MaxMom for elements with an aspect
ratio over MaxAR prevents the generation of badly shaped elements
parallel to the traces. The MaxPnt strategy produces a mesh quite similar
to the mesh produced by MaxMom due to the fact that we fave few cells
with a number of vertices larger than MaxNP during the refinement process.
In Figures
\ref{fig:ten1f72sAr} and \ref{fig:ten4f72Ar} we report the minimum,
the mean and the maximum aspect ratios of the cells on Fracture 72 along the
refinement process. We remark that in all the strategies we use MaxMom
strategy to refine the elements with large aspect ratio. The MaxMom
and the MaxEdg strategies produce a decreasing mean aspect ratio,
whereas TrDir and MaxPnt have a different behaviour. A slight
difference from the two figures can be seen in the TrDir plots. The
average AR grows in the DFN86E04 test case because the
MaxMom strategy is less used due to the weaker refinement around the traces
due to the different transmissivities on the intersecting fractures
that justify a smaller flux (notice the more coarse mesh for DFN86E04 comparing Figures \ref{fig:ten1f72slast} and \ref{fig:ten4f72slast}).



\bibliographystyle{elsarticle-harv}
\bibliography{scico2mp,dfn,articolo,VEMbibl}

\end{document}